\documentclass[11pt]{article}
\usepackage{graphicx}
\usepackage{latexsym}
\usepackage{amssymb}

\def\beq{\begin{eqnarray}}
\def\eeq{\end{eqnarray}}
\def\I{\mbox{\bf I}}
\def\mod{\mbox{mod}\ }
\def\met{_{\mbox{\tiny met}}}
\def\seq{_{\mbox{\tiny seq}}}

\begin{document}

\fontsize{11}{14.5pt}\selectfont

\begin{center}

{\small Technical Report No.\ 0506,
 Department of Statistics, University of Toronto}

\vspace*{0.9in}

{\LARGE \bf The Short-Cut Metropolis Method}\\[16pt]

{\large Radford M. Neal}\\[3pt]
 Department of Statistics and Department of Computer Science \\
 University of Toronto, Toronto, Ontario, Canada \\
 \texttt{http://www.cs.utoronto.ca/$\sim$radford/} \\
 \texttt{radford@stat.utoronto.ca}\\[10pt]

 2 August 2005

\end{center}

\vspace*{8pt}

\noindent \textbf{Abstract.}  I show how one can modify the
random-walk Metropolis MCMC method in such a way that a sequence of
modified Metropolis updates takes little computation time when the
rejection rate is outside a desired interval.  This allows one to
effectively adapt the scale of the Metropolis proposal distribution,
by performing several such ``short-cut'' Metropolis sequences with
varying proposal stepsizes.  Unlike other adaptive Metropolis schemes,
this method converges to the correct distribution in the same fashion
as the standard Metropolis method.

\section{\hspace*{-7pt}Introduction}\label{sec-intro}\vspace*{-10pt}

The Metropolis algorithm of Metropolis, Rosenbluth, Rosenbluth,
Teller, and Teller (1953) is the first and perhaps still the most
widely-used Markov chain Monte Carlo (MCMC) method.  The Metropolis
procedure defines a probabilistic transition from a state $x$ (usually
of high dimension) to a state $x'$ that leaves some desired
distribution, with density function $\pi(x)$, invariant.  Starting
from any initial state, a Markov chain using such transitions will
produce states that asymptotically converge in distribution to $\pi$,
provided that the transitions are capable of moving from any region of
the state space to any other.  States from the latter parts of one or
more simulations of such a Markov chain (called ``runs'') are then
used to make Monte Carlo estimates for expectations of functions of
state with respect to $\pi$.  For more information on Metropolis and
other MCMC methods, and their applications in statistics and
statistical physics, see (for example) Neal (1993), Tierney (1994),
and Liu (2001).

The Metropolis algorithm requires that we specify a suitable
``proposal'' distribution, which may depend on the current state,
whose conditional density, $g(x^*|x)$, satisfies the symmetry
condition $g(x^*|x)=g(x|x^*)$.  A transition (or ``update'') from
state $x$ is performed by generating an $x^*$ according to
$g(\cdot|x)$, and then ``accepting'' $x^*$ as the new state, $x'$,
with probability
\beq
   a(x,x^*) & = & \min\big[\, 1,\ \pi(x^*)/\pi(x)\,\big]
\label{eq-met-acc}\eeq
If $x^*$ is not accepted (ie, is ``rejected''), we let the new state, 
$x'$, be the same as the old state, $x$.  

The most generally-useful class of proposal distributions take the
form of adding a random offset to the current state --- ie, $x^* =
x+w\delta$.  Here, $w$ is a scalar ``stepsize'' parameter, and
$\delta$ is a vector drawn from some distribution not depending on
$x$, with density function $f(\delta)$, which must be symmetrical
around zero (ie, $f(\delta)=f(-\delta)$).  In this paper, I will
mostly consider distributions for $\delta$ in which all components are
non-zero, so that the proposed $x^*$ differs in all components from
the current $x$.  Metropolis procedures in which only one (or a small
group) of components are changed in each update are also common, and
will be discussed briefly at the end.

Such ``random walk'' Metropolis updates can explore complex
distributions whose shape is not known \textit{a~priori}.  However,
the efficiency of this exploration depends critically on a proper
choice for the stepsize, $w$.  If $w$ is too large, we may find that
almost all proposals are rejected, leading to very inefficient
exploration.  But if $w$ is too small, each update will change the
state by only a small amount, so many updates will be required to move
a substantial distance (an effect exacerbated by the fact that
movement is via a random walk).

We often do not know enough about $\pi$ initially to choose an
appropriate value for the stepsize~$w$.  A common approach is to
perform some preliminary runs using various values for $w$, and then
use statistics from these runs to select a value for $w$ that appears
to give a rejection rate between 50\% and 80\%.  Alternatively,
we can manually change $w$ during the early part of a run until we
believe we've found a suitable value.  With either method, we then use
the chosen value for $w$ in the remainder of the run, or in one or
more new runs, states from which are used for estimating expectations.
Although many useful results have been obtained in this way, this
procedure seems a bit wasteful and inelegant.  It also may not work
well in situations where different stepsizes are appropriate in
different regions of the state space.  Accordingly, many people have
sought methods in which the stepsize is automatically adapted in some
suitable way throughout the course of a run.

One obvious approach is to continually change the stepsize based on
the rejection rate in past updates.  Unfortunately, many naive methods
of this sort will give wrong answers.  For example, we might for each
update use one of two values for the stepsize, $w_0$ or $w_1$ (with
$w_0<w_1$), choosing $w_0$ if the rejection rate in the last 10
updates was greater than 50\%, and $w_1$ otherwise.  There is no
reason to think that this will produce the right answer, however,
since the dependence on past updates destroys the Markov property of
the process, which undermines the proof that the distribution of the
state converges to $\pi$.  Indeed, this and similar methods will tend
to spend too much time in regions of the state space where a small
stepsize is needed to achieve a small-enough rejection rate, since
when the Markov chain is in such a region, it will effectively be
operating at a slower pace than in regions where a larger stepsize is
chosen.

Estimates that are asymptotically correct can be obtained if changes
to the stepsize become smaller as the run progresses.  Methods of this
sort have been investigated by Haario, Saksman, and Tamminen (2001),
Atchad\'{e} and Rosenthal (2005), and Andrieu and Moulines (2005), who
show that when using certain adaptation schemes, estimates of
expectations for functions of state will converge to the correct
values under certain conditions.  These conditions may not be easy to
satisfy or to verify, however, and even if they can be shown to hold,
properly assessing the error in the estimates found will be more
complex than for ordinary MCMC estimation, for which assessing
accuracy is already a non-trivial problem, generally requiring some
human judgement.  More fundamentally, since the correctness of these
methods depends on the stepsize becoming more and more stable as the
run progresses, it is unclear what advantage these methods might have
over the much simpler procedure of adapting the stepsize only during
an initial portion of the run, fixing it after some number of updates
--- essentially an automated version of one of the commonly-used
manual procedures described above.  Choosing a suitable point at which
to fix the stepsize (either manually or automatically) seems easier
than assessing the degree to which the accuracy of estimates may have
been affected by continual changes in the stepsize.

We would in any case prefer a method that can use different
stepsizes in different regions of the state space.  Something like
this is the objective of the ``delayed rejection'' Metropolis method
(Tierney and Mira 1999; Green and Mira 2001).  In this method, an
update may consider several proposals before one is accepted, and
later proposals may be influenced by the results of earlier proposals
within the same update.  ``Slice sampling'' updates (Neal 2003)
achieve a similar effect, in a manner that may be more useful in
practice.  However, these methods adapt only temporarily, as part of a
larger update.  Adaptation must begin over again at the start of the next
update.

In this paper, I describe a new approach to ``adapting'' the stepsize,
which is best understood by first considering the following wasteful
strategy.  Suppose we are uncertain whether to use $w=1$ or $w=10$.
We might therefore use both of these stepsizes in turn.  For instance,
we might alternate between performing 30 updates with $w=1$ and
performing 90 updates with $w=10$.  Although this procedure may avoid
disaster by ensuring that we use a good stepsize at least some of the
time, it will be reasonably efficient only when $w=10$ is the best
stepsize, in which case we will be using the right stepsize for
$90/(90+30) = 3/4$ of the time.  When $w=1$ is the best step size, but
$w=10$ is too large, almost all the 90 updates done with $w=10$ will
be rejected, so only $30/(90+30)=1/4$ of the updates will be useful,
which seems less than satisfactory.

Suppose, however, that we could somehow arrange that simulating a
sequence of Metropolis updates done using a stepsize that is too
large, and which therefore are almost always rejected, takes little
computation time.  In particular, suppose that simulating a sequence
of $K$ Metropolis updates with a stepsize that almost always leads to
rejection can be done in the same time as 10 ordinary Metropolis
updates, \textit{regardless} of how large $K$ is.  With this
``short-cut'' method, the strategy described above is attractive even
when $w=1$ is the best stepsize --- the 90 updates with $w=10$ are
done in the same time as 10 ordinary updates, so $30/(10+30)=3/4$ of
the time is devoted to updates with the best stepsize.

A variation on this approach is possible if a we can also take a
short-cut when simulating a sequence of Metropolis updates in which
the rejection rate is much smaller than desired.  Suppose that a
sequence of $K$ updates can be effectively reduced to only 10 updates
when almost no updates are rejected, or when almost all updates are
rejected, regardless of how large $K$ is.  We could then alternately
perform a sequence of 60 updates with $w=1$ and a sequence of 60
updates with $w=10$.  If only one of these values for $w$ produces a
reasonable rejection rate (not near 0 or 1), then 60/(10+60)=6/7 of
our time will be spent simulating the useful updates, with only 1/7 of
our time wasted on updates for which the stepsize is too large or too
small.

From a computational point of view, these short-cut Metropolis
updates can be used to adaptively change the stepsize used during the
bulk of the computation, effectively allowing different stepsizes to
be used in different regions of the state space.  From a mathematical
point of view, however, we are still simply performing sequences of
updates using pre-determined stepsizes.  The Markov property still
holds, and all the usual MCMC theorems regarding convergence of the
distribution to $\pi$ and of estimates for expectations with respect
to $\pi$ to their correct values still apply.

The reader may wonder, of course, how one might manage to simulate $K$
Metropolis updates in time that does not depend on $K$.  To see how
the Metropolis algorithm can be modified to achieve this, we first
need to re-interpret a Metropolis update as a deterministic
transformation.  This is the topic of the next section.

\section{\hspace*{-7pt}A deterministic view of the Metropolis 
         algorithm}\label{sec-det}\vspace*{-10pt}

In this section I will describe how a standard Metropolis update can
be viewed as a deterministic transformation following a stochastic
extension of the state space to include auxiliary variables.  This
sets the stage for the modification in the next section, which allows
a sequence of Metropolis updates with a badly-chosen stepsize to be
simulated very quickly.  Note that the complexities of viewing
Metropolis updates in this way need not be reflected in the
implementation --- in this section, the actual implementation is just
the standard Metropolis algorithm.  Viewing Metropolis updates in
terms of auxiliary variables and deterministic transformations is
necessary only to show that the modified updates still leave the
desired distribution invariant.

To view a Metropolis update as a deterministic transformation, we need
to temporarily introduce auxiliary variables, a technique that is
familiar from other MCMC methods as well, such as slice sampling (Neal
2003).  Suppose we wish to sample $x$ from distribution $\pi(x)$,
using a chain which has this as its unique invariant distribution.
Rather than directly defining an update of $x$ that leaves $\pi(x)$
invariant, we can expand the state space to pairs $(x,y)$, define a
distribution $\pi(x,y) = \pi(x)\pi(y|x)$, where $\pi(x)$ is the
original distribution of interest, and find a way of updating the pair
$(x,y)$ that leaves $\pi(x,y)$ invariant.  We can introduce the
auxiliary variable $y$ temporarily by sampling a value for it from
$\pi(y|x)$, performing an update of $(x,y)$, and then discarding $y$.
It is easy to see that this procedure leaves $\pi(x)$ invariant if the
update for $(x,y)$ leaves $\pi(x,y)$ invariant.

Deterministic updates can form part of a valid MCMC scheme as long as
they leave the desired distribution invariant.  Suppose the desired
distribution for the state, $z$, has density $\pi(z)$.  Let $T(z)$ be
a 1-to-1 mapping from the state space onto itself.  If $z$ has density
$\pi(z)$, and $z'=T(z)$, then according to the standard formula for
transformation of probability densities, the density of $z'$ is
$\pi(T^{-1}(z'))\,/\,|\det T'(T^{-1}(z'))\,|$, where $T'(z)$ is the
Jacobian matrix for the transformation.  It follows that the
transformation $z'=T(z)$ will leave $\pi$ invariant as long as
$\pi(T(z))=\pi(z)$ and $|\det T'(z)\,|=1$ for all $z$.  If $z$ is
partially discrete, the Jacobian condition applies only to the
continuous portion.

Consider a standard random-walk Metropolis update of $x$, with
invariant density $\pi(x)$, symmetrical proposal density $f(\delta)$,
and stepsize $w$.  To view this update in terms of auxiliary variables
and a deterministic transformation, we first expand the state space by
introducing two auxiliary variables:\ \ $\delta$, a vector of the same
dimension as $x$, and $e$, a positive real number.  We
then define the joint density $\pi(x,\delta,e) =
\pi(x)\,f(\delta)\,\exp(-e)$, in which $x$, $\delta$, and $e$ are independent.
Finally, we define a transformation, $T\met$, from 
$(x,\,\delta,\,e)$ to $(x',\delta',e')$ as follows:\vspace*{4pt}
\beq
  T\met(x,\,\delta,\,e) & = & \left\{\begin{array}{ll}
    (x\!+\!w\delta,\ -\delta,\ e + \log(\pi(x\!+\!w\delta)/\pi(x))\ 
       & \mbox{if $e + \log(\pi(x\!+\!w\delta)/\pi(x)) > 0$}\ \ \\[3pt]
    (x,\,\delta,\,e)          
       & \mbox{otherwise}
  \end{array}\right.\label{eq-met-trans}\\[-7pt]\nonumber
\eeq
One can easily see that this transformation is 1-to-1 and onto, since
for any point $(x',\delta',e')$, exactly one point will transform to
this point, namely the point $(x'\!+\!w\delta',-\delta',
e' - \log(\pi(x')/\pi(x'\!+\!w\delta')) )$ if 
$e' - \log(\pi(x')/\pi(x'\!+\!w\delta'))>0$ or the point
$(x',\delta',e')$ otherwise.  The Jacobian matrix for this transformation
is the identity if $e + \log(\pi(x\!+\!w\delta)/\pi(x)) \le 0$, and is
otherwise
\beq
  \left[\begin{array}{ccc}
     \I            & w\I          & 0 \\[3pt]
      \mbox{\bf 0} & -\I          & 0 \\[3pt]
      \mbox{\bf a} & \mbox{\bf ~b} & 1
  \end{array}\right]
\eeq
where $\I$ is the identity matrix with dimensions equal to the number of
components in $x$ and $\delta$, and the values of \textbf{a} and 
\textbf{b} are irrelevant.  One can easily see that the absolute value
of the determinant of this matrix is one, since the only non-zero term
in the determinant is the product of the entries on the main diagonal,
which are all $1$ or $-1$.  Finally, the probability density of 
$(x',\delta',e')$ is the same as that of $(x,\delta,e)$, since when
these points differ
\beq
\pi(x',\delta',e') & = & 
  \pi(x\!+\!w\delta,\ -\delta,\ e + \log(\pi(x\!+\!w\delta)/\pi(x)) \\[5pt]
& = &
  \pi(x\!+\!w\delta)\, f(-\delta)\, \exp(-(e + \log(\pi(x\!+\!w\delta)/\pi(x))) 
  \\[5pt]
& = & 
  \pi(x)\, f(\delta)\, \exp(-e) \ \ =\ \ \pi(x,\delta,e)
\eeq

We can therefore conclude that the following procedure for updating 
$x$ to $x'$ leaves $\pi(x)$ invariant:\vspace{12pt}\\
\hspace*{2.6in}\textbf{Procedure 1:}\vspace*{-7pt}
\begin{enumerate}
\item[1)] Randomly pick values for the temporary auxiliary variables,
          as follows:\vspace*{-4pt}
 \begin{enumerate}
 \item[a)] Pick a value for $\delta$ according to the density function 
           $f(\delta)$.
 \item[b)] Pick a value for $e$ from the exponential distribution with
           mean one, whose density function is $\exp(-e)$.
 \end{enumerate}
\item[2)] Set $(x',\,\delta',\,e')\ =\ T\met(x,\,\delta,\,e)$, where $T\met$ 
          is defined by equation (\ref{eq-met-trans}).
\item[3)] Forget $\delta'$ and $e'$, leaving just $x'$ as the new 
          state.\vspace*{-6pt}
\end{enumerate}
The distribution sampled in Step (1) is $\pi(\delta,e|x)$, so the joint
distribution at this point will be $\pi(x,\delta,e)$ if the distribution
of $x$ was $\pi(x)$ before.  Step (2) leaves invariant this joint distribution,
and hence also the marginal distribution for $x$, which we return to in 
Step (3).

The effect of this procedure is equivalent to a standard random-walk
Metropolis update.  In step~(3), $x'$ will be set to either $x$ or
$x\!+\!w\delta$. The probability of the latter is the probability that
a value drawn from an exponential distribution with mean one is
greater than $-\log(\pi(x\!+\!w\delta)/\pi(x))$, which is
$\min[1,\,\pi(x\!+\!w\delta)/\pi(x)]$, the same as the Metropolis
acceptance probability of (\ref{eq-met-acc}).

We can also view an entire sequence of $K$ random-walk Metropolis updates 
in terms of an initial stochastic extension of the state space
to include auxiliary variables followed by a single deterministic
transformation.  We now need $K$ copies of the auxiliary variables
used above, denoted by $\delta_0,\ldots,\delta_{K-1}$ and
$e_0,\ldots,e_{K\!-\!1}$, plus an index variable, $i$, in 
$\{0,\ldots,K\!-\!1\}$.  The joint density for $x$ and these auxiliary 
variables is defined to be\vspace*{-4pt}
\beq
  \pi(x,i,\delta_0,\ldots,\delta_{K-1},e_0,\ldots,e_{K-1}) & = & 
  \pi(x)\, {1 \over K}\, 
  \prod_{j=0}^{K-1} f(\delta_j) \,\prod_{j=0}^{K-1} \exp(-e_j) 
\eeq
In other words, the auxilary variables are all independent of $x$ and
each other, and their distributions are as 
previously defined, except for $i$, which is given a uniform distribution.
We will repeatedly apply a transformation in which
$(x,i,\delta_0,\ldots,\delta_{K-1},e_0,\ldots,e_{K-1})$ is mapped to
the point found by transforming $(x,\delta_i,e_i)$ to $(x',\delta_i',e_i')$
according to the mapping $T\met$ of equation~(\ref{eq-met-trans}), the
value of $i$ is change to $i'=i+1\ (\mod K)$, 
and the other auxiliary variables are left unchanged.   As shown above,
the $T\met$ part of this transformation is 1-to-1 
and onto, the absolute value of the determinant of its Jacobian matrix 
is one, and the probability density of the transformed point is the same 
as that of the original point.  The $i'=i+1\ (\mod K)$ part is also 1-to-1 
and onto, and since $i$ is discrete, there is no Jacobian to worry about.  
Since $i$ has a uniform distribution, the density is also not changed by this
part of the transformation.  

The entire transformation therefore leaves $\pi$ invariant.  We can
write the procedure using this transformation in detail as 
follows:\vspace{12pt}\\
\hspace*{2.6in}\textbf{Procedure 2:}\vspace*{-7pt}
\begin{enumerate}
\item[1)] Randomly pick values for the temporary auxiliary variables,
          as follows:\vspace*{-4pt}
 \begin{enumerate}
 \item[a)] Pick values for $\delta_0,\ldots,\delta_{K-1}$ independently 
           according to the density function $f(\delta)$.
 \item[b)] Pick values for $e_0,\ldots,e_{K-1}$ independently from the 
           exponential distribution with
           mean one, whose density function is $\exp(-e)$.
 \item[c)] Pick a value for $i$ uniformly from 
           $\{0,\ldots,K\!-\!1\}$.
 \end{enumerate}
\item[2)] Transform $x$ and the auxiliary variables by repeating the 
          following sequence of transformations $K$ times:\vspace*{-4pt}
 \begin{enumerate}
 \item[a)] Apply the transformation $T\met$ to $(x,\delta_i,e_i)$
 \item[b)] Add one to $i$, modulo $K$.
 \end{enumerate}
\item[3)] Forget $\delta_0,\ldots,\delta_{K-1}$, $e_0,\ldots,e_{K-1}$, and
          $i$, leaving just $x'$ as the new state.\vspace*{-6pt}
\end{enumerate}
It easy to see that this procedure is equivalent to simply performing 
$K$ successive random-walk Metropolis updates.  The index $i$ takes on
each of its $K$ possible values exactly once, and for each value of 
this index, a Metropolis update is performed utilizing the auxiliary
variables $\delta_i$ and $e_i$.

\section{\hspace*{-7pt}Short-cut simulation of sequences of Metropolis 
                       updates}\label{sec-mod}\vspace*{-10pt}

Consider now a sequence of $K$ Metropolis updates, divided into $M$
groups of $L$ updates (so that $K=ML$).  I show here how we can look
at the number of rejections within each group as we simulate it, and
modify subsequent actions so that no further computation is required
once two groups are simulated in which the number of rejections is
outside some desired range.  For example, if we set $L=5$, we can
achieve the result mentioned in the introduction --- if the stepsize
is so large that almost all Metropolis updates are rejected, or so
small that almost none are rejected, the $K$ updates can usually be
computed using the time normally required for only 10 updates,
regardless of how large $K$ is.

First, however, let's once again look at how we can perform standard
Metropolis updates using auxiliary variables and a deterministic
transformation, this time expressed in terms of $M$ successive
transformations, each of which corresponds to $L$ Metropolis updates.  
We introduce a new auxiliary variable, $s$, 
in $\{-1,+1\}$, which will determine in which direction the index $i$ moves.  
It will be uniformly distributed,
independently of the other variables, so the joint distribution of
$x$ and all auxiliary variables is now
\beq
  \pi(x,i,s,\delta_0,\ldots,\delta_{K-1},e_0,\ldots,e_{K-1}) & = & 
  \pi(x)\, {1 \over K}\, {1 \over 2}\,
  \prod_{j=0}^{K-1} f(\delta_j) \,\prod_{j=0}^{K-1} \exp(-e_j) 
\eeq

We can use arguments paralleling those used in the previous section 
to justify Procedure~2 to show that the following procedure for updating 
$x$ leaves $\pi$ invariant:\vspace{12pt}\\
\hspace*{2.6in}\textbf{Procedure 3:}\vspace*{-7pt}
\begin{enumerate}
\item[1)] Randomly pick values for the temporary auxiliary variables,
          as follows:\vspace*{-4pt}
 \begin{enumerate}
 \item[a)] Pick values for $\delta_0,\ldots,\delta_{K-1}$ independently 
           according to the density function $f(\delta)$.
 \item[b)] Pick values for $e_0,\ldots,e_{K-1}$ independently from the 
           exponential distribution with
           mean one, whose density function is $\exp(-e)$.
 \item[c)] Pick a value for $i$ uniformly from 
           $\{0,\ldots,K\!-\!1\}$.
 \item[d)] Pick a value for $s$ uniformly from $\{-1,+1\}$.
 \end{enumerate}
\item[2)] Transform $x$ and the auxiliary variables by repeating the 
          following sequence of transformations $M$ times:\vspace*{-4pt}
 \begin{enumerate}
 \item[a)] First, apply the transformation $T\met$ to $(x,\delta_i,e_i)$.
           Second, do the following $L-1$ times:\ \ add $s$ to $i$, 
           modulo $K$, and then apply $T\met$ to $(x,\delta_i,e_i)$.  
           Third, negate~$s$.
 \item[b)] Negate $s$.
 \item[c)] Add $s$ to $i$, modulo $K$.
 \end{enumerate}
\item[3)] Forget $\delta_0,\ldots,\delta_{K-1}$, $e_0,\ldots,e_{K-1}$, 
          $i$, and $s$, leaving just $x'$ as the new state.\vspace*{-6pt}
\end{enumerate}
Note that each of Steps (2a), (2b), and (2c) leaves $\pi$ invariant.
The negation of $s$ in Step~(2b) simply undoes the negation at the 
end of Step~(2a), so the only difference between this procedure and the
one in the preceding section is that the initial value of $s$ determines
whether the $\delta_i$ and $e_i$ are utilized in increasing order 
or decreasing order.  The order makes no difference, so this procedure is 
also equivalent to performing $K=ML$ standard random-walk Metropolis
updates.

We can now modify this procedure to avoid computation when the
rejection rate is greater or less than desired.  The modification will
have the effect of negating $s$, thereby reversing the direction in
which $i$ moves, whenever the number of rejections within a group of
$L$ updates is outside some desired range.  The updates following this
will undo earlier updates, leading back to previously computed states,
which need not be computed again.  Once the original state is reached,
new states will be computed, but if a second reversal occurs, all
subsequent updates will result in states that have already been
computed.  An extreme, but not necessarily uncommon, case occurs when
the number of rejections in the first group of $L$ updates is outside
the desired range, and the same is true of the group of $L$ updates
simulated next.  The final state will then be the same as the original
state.  A total of only $2L$ updates requiring computation will have
been done, regardless of $K$.

To see how this modification can be done, first note that the
transformation in Step~(2a) --- call it $T\seq$ --- is its own
inverse.  $T\seq$ will change $i$ to $i' = i+(L\!-\!1)s\ (\mod K)$ and
$s$ to $s'= -s$.  Applying $T\seq$ a second time would change $i$ to
$i'' = i'+(L\!-\!1)s'\ (\mod K) = i$ and change $s$ to $s''=-s'=s$.
Changes to $x$ and the $\delta_i$ and $e_i$ would also be undone.  The
values of $i$ used in these updates would be visited in reverse order.
If the original update was a rejection, for which
$T\met(x,\delta_i,e_i) = (x,\delta_i,e_i)$, it will be a rejection in
this reverse pass as well, leaving the original state unchanged.  If
the original update instead had $T\met(x,\delta_i,e_i) =
(x\!+\!w\delta_i,\ -\delta_i,\ e_i +
\log(\pi(x\!+\!w\delta_i)/\pi(x))$, applying $T\met$ a second time
will again restore the original state of $(x,\delta_i,e_i)$.

In general, any transformation, $z' = T(z)$, that is its own inverse
(and hence is also 1-to-1 and onto), that has a Jacobian matrix for
which the absolute value of the determinant is one, and that satisfies
\mbox{$\pi(T(z))=\pi(z)$} for all $z$ will, as we have seen, leave
the density $\pi(z)$ invariant.  For any set $A$, consider the modified
transformation, $\tilde T^A$, defined by \beq
  \tilde T^A(z) & = & \left\{\begin{array}{ll}
     z & \mbox{if $z\in A$ or $T(z) \in A$} \\[3pt]
     T(z) & \mbox{otherwise}
  \end{array}\right.
\label{eq-transA}\eeq
One can easily show that $\tilde T^A$ is also its own inverse:  If $z \in A$
or $T(z) \in A$, then $\tilde T^A(\tilde T^A(z)) = \tilde T^A(z) = z$,
and if $z \notin A$ and $T(z) \notin A$, then $T(T(z))=z \notin A$,
and $\tilde T^A(\tilde T^A(z)) = \tilde T^A(T(z)) = T(T(z)) = z$.
One can also easily see that the determinant of the Jacobian matrix of
$\tilde T^A$ has absolute value one, and that $\pi(\tilde
T^A(z))=\pi(z)$.  It follows that $\tilde T^A$ leaves $\pi$ invariant.

If we modify Step~(2a) in Procedure~3 to use $T\seq^A$ rather
than $T\seq$, the overall procedure will therefore still leave $\pi$
invariant.  A variety of choices for the set $A$ may be useful, but
for the moment I will consider only sets, $R_{l,h}$, consisting of
values for $x$ and the auxiliary variables for which applying $T\seq$
will result in $L$ updates in which the number of rejections is
outside the interval $[l,h]$.  We might, for instance, use the set
$R_{1,L-1}$, in order to pick out groups in which either all updates
are rejected or no updates are rejected, or the set $R_{0,L-1}$, in
order to pick out only groups in which all updates are rejected.  By
considering the details of Step~(2a), one can see that a point is in
$R_{l,h}$ if and only if the image of this point under $T\seq$ is also in
$R_{l,j}$, since if an update is rejected when simulating in the forward
direction, it will also be rejected when simulating backwards.

Modifying Procedure~3 to use $T\seq^{R_{l,h}}$ rather than $T\seq$ in Step~(2a),
and then merging this step with Step~(2b), yields the following 
procedure:\vspace{12pt}\\
\hspace*{2.0in}\textbf{Procedure 4 (Short-Cut Metropolis):}\vspace*{-7pt}
\begin{enumerate}
\item[1)] Randomly pick values for the temporary auxiliary variables,
          as follows:\vspace*{-4pt}
 \begin{enumerate}
 \item[a)] Pick values for $\delta_0,\ldots,\delta_{K-1}$ independently 
           according to the density function $f(\delta)$.
 \item[b)] Pick values for $e_0,\ldots,e_{K-1}$ independently from the 
           exponential distribution with
           mean one, whose density function is $\exp(-e)$.
 \item[c)] Pick a value for $i$ uniformly from 
           $\{0,\ldots,K\!-\!1\}$.
 \item[d)] Pick a value for $s$ uniformly from $\{-1,+1\}$.
 \end{enumerate}
\item[2)] Transform $x$ and the auxiliary variables by repeating the 
          following sequence of transformations $M$ times:\vspace*{-4pt}
 \begin{enumerate}
 \item[a)] First, apply the transformation $T\met$ to $(x,\delta_i,e_i)$.
           Second, do the following $L-1$ times:\ \ add $s$ to $i$, 
           modulo $K$, and then apply $T\met$ to $(x,\delta_i,e_i)$.  
           Finally, if the number of these applications of $T\met$ 
           that were rejections is less than $l$ or greater than $h$,
           change $i$ back to its value at the start of this step (ie, subtract
           $(L\!-\!1)s$ from it, modulo $K$), change $x$ back to its
           value at the start of this step, and negate $s$.
 \item[b)] Add $s$ to $i$, modulo $K$.
 \end{enumerate}
\item[3)] Forget $\delta_0,\ldots,\delta_{K-1}$, $e_0,\ldots,e_{K-1}$, 
          $i$, and $s$, leaving just $x'$ as the new state.\vspace*{-6pt}
\end{enumerate}
The new Step (2a) above reverses the direction, $s$, only when a group
is simulated in which the number of rejections is outside the desired
range.  Procedure~4 will therefore mimic a sequence of standard
Metropolis updates as long as the number of rejections in each group
of $L$ updates is in this range.  If a group of updates fails this test,
the direction reverses, and previous states are revisited, without the
need to recompute them.  Once these previous states have all been
visited, new states are again simulated.  If a second group of updates
fails the test, the direction reverses again, and states are again
revisited, with no further computation of states being required.  It
is possible that the direction will subsequently reverse many times,
with states being revisited again and again in a back-and-forth
manner.  Note that if the numbers of rejections in the first two
groups of updates are outside the desired range, we simply stay at 
the initial state.

\begin{figure}[p]

\includegraphics[width=6.4in]{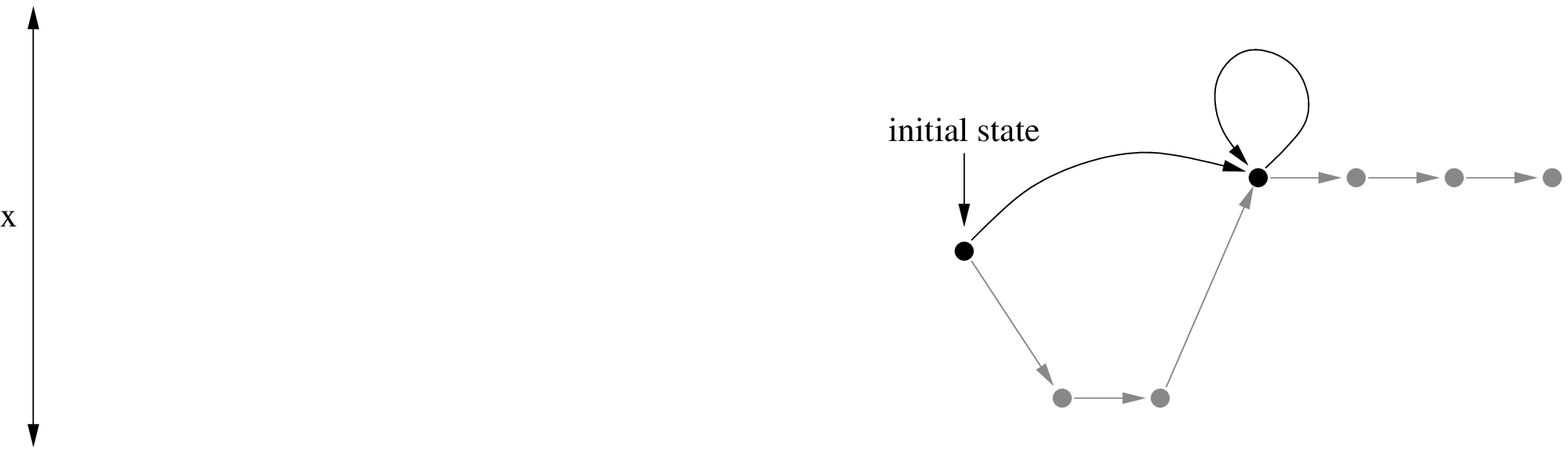}

\hspace*{24pt}%
\makebox[29.3pt]{~}%
\makebox[29.3pt]{~}%
\makebox[29.3pt]{~}%
\makebox[29.3pt]{~}%
\makebox[29.3pt]{~}%
\makebox[29.3pt]{~}%
\makebox[29.3pt]{~}%
\makebox[29.3pt]{~}%
\makebox[29.3pt]{~}%
\makebox[29.3pt]{$\scriptstyle \delta_9$}%
\makebox[29.3pt]{$\scriptstyle \delta_{10}$}%
\makebox[29.3pt]{$\scriptstyle \delta_{11}$}%
\makebox[29.3pt]{$\scriptstyle \delta_{12}$}%
\makebox[29.3pt]{$\scriptstyle \delta_{13}$}%
\makebox[29.3pt]{$\scriptstyle \delta_{14}$}%

\vspace*{-3pt}
\rule{6.4in}{0.2pt}
\vspace*{18pt}

\includegraphics[width=6.4in]{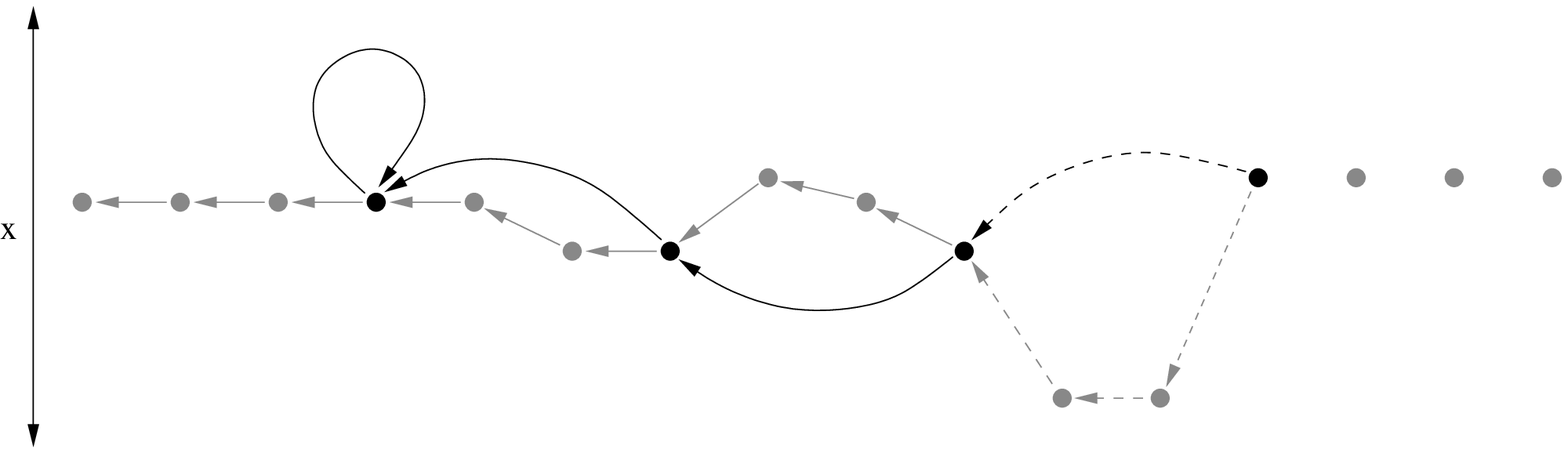}

\hspace*{24pt}%
\makebox[29.3pt]{$\scriptstyle \delta_0$}%
\makebox[29.3pt]{$\scriptstyle \delta_1$}%
\makebox[29.3pt]{$\scriptstyle \delta_2$}%
\makebox[29.3pt]{$\scriptstyle \delta_3$}%
\makebox[29.3pt]{$\scriptstyle \delta_4$}%
\makebox[29.3pt]{$\scriptstyle \delta_5$}%
\makebox[29.3pt]{$\scriptstyle \delta_6$}%
\makebox[29.3pt]{$\scriptstyle \delta_7$}%
\makebox[29.3pt]{$\scriptstyle \delta_8$}%
\makebox[29.3pt]{$\scriptstyle -\delta_9$}%
\makebox[29.3pt]{$\scriptstyle -\delta_{10}$}%
\makebox[29.3pt]{$\scriptstyle -\delta_{11}$}%
\makebox[29.3pt]{~}%
\makebox[29.3pt]{~}%
\makebox[29.3pt]{~}%

\vspace*{-3pt}
\rule{6.4in}{0.2pt}
\vspace*{18pt}

\includegraphics[width=6.4in]{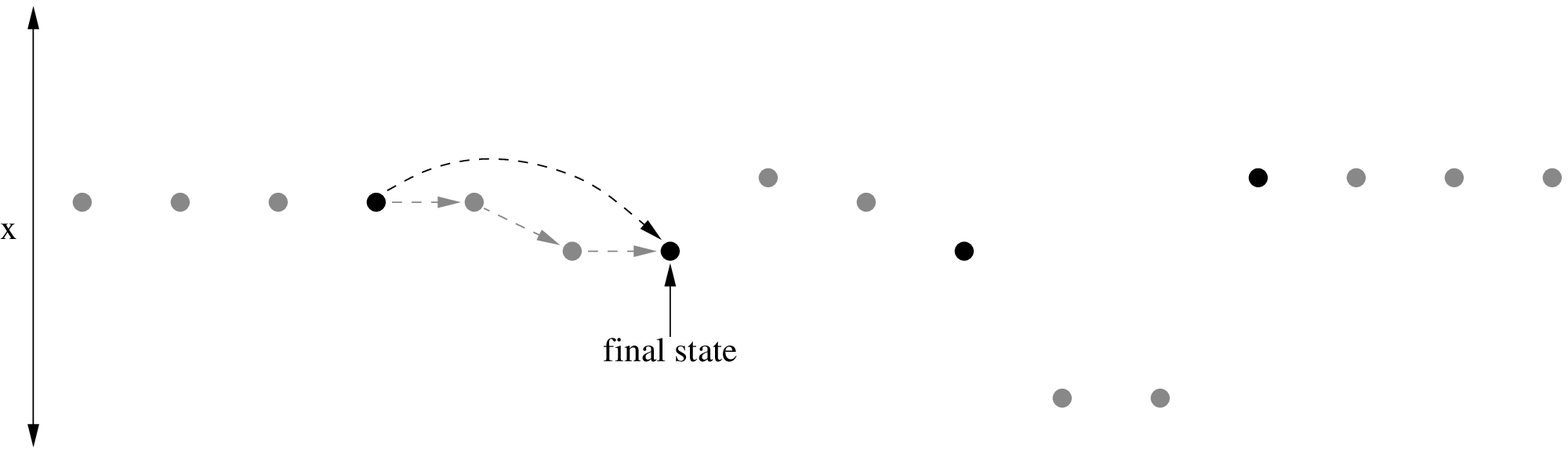}

\hspace*{24pt}%
\makebox[29.3pt]{~}%
\makebox[29.3pt]{~}%
\makebox[29.3pt]{~}%
\makebox[29.3pt]{$\scriptstyle -\delta_3$}%
\makebox[29.3pt]{$\scriptstyle -\delta_4$}%
\makebox[29.3pt]{$\scriptstyle -\delta_5$}%
\makebox[29.3pt]{~}%
\makebox[29.3pt]{~}%
\makebox[29.3pt]{~}%
\makebox[29.3pt]{~}%
\makebox[29.3pt]{~}%
\makebox[29.3pt]{~}%
\makebox[29.3pt]{~}%
\makebox[29.3pt]{~}%
\makebox[29.3pt]{~}%

\vspace*{2pt}

\caption[]{An illustration of the short-cut Metropolis method, with $L=3$,
$M=7$, $l=0$, and $h=2$.  The vertical
axis represents the state (here one dimensional).  The horizontal axis 
represents the index, $i$, of the auxiliary variables, $(\delta_i,e_i)$,
used for each update.  The value of $\delta$ for each 
update is shown below --- either the original $\delta_i$ or its negation. 
Grey arrows show the $K$ applications of $T\met$; dotted arrows
indicate applications that produce already-computed states.
Black arrows show the results of the full transformation defined by 
Step~(2a) of Procedure~4, with dotted arrows again 
indicating the results have already been computed.  The top panel 
shows the first two groups of updates, with the second group consisting 
only of rejections.  The middle panel shows the next group revisiting
the first group of updates in reverse order, and then three additional 
groups, simulated starting with the original state, with the last group 
consisting only of rejections.  The bottom panel shows the final group 
of $L$ updates, which again revisits already computed states.}\label{fig-ill}

\end{figure}

Figure~\ref{fig-ill} illustrates this short-cut Metropolis procedure,
for the case where $L=3$ and $M=7$, and hence $K=ML=21$, and where
$l=0$ and $h=L\!-\!1$ --- ie, we
reverse direction when all updates in a group are rejected, but not
when none of the updates are rejected.

\section{\hspace*{-7pt}Implementation issues}\label{sec-imp}\vspace*{-10pt}

When implementing the short-cut Metropolis method, the description in
Procedure~4 can be simplified.  There is no need to randomly pick a
value for $i$ --- this step is needed for the proof of validity, but
one can see that the distribution of the final result is the same
regardless of which value for $i$ was picked.  Time can be saved by
generating values for $\delta_i$ and $e_i$ only if and when they are
needed.  Once used, the values for these auxiliary variables can be
forgotten, provided the states that were generated using them are
saved.  

Most importantly, as was discussed above, although Procedure~4 as written
appears to involve $K \!=\! ML$ applications of $T\met$ (defined by
equation~(\ref{eq-met-trans})), each requiring evaluation of $\pi(x)$
at two states, when the direction of simulation reverses due to a group
of $L$ rejections, many of these applications of $T\met$ can be
avoided, since they produce states that have already been computed.
Furthermore, if the value of $\pi(x)$ is saved for all states that
have been computed, an application of $T\met$ will require only the
evaluation of $\pi(x+w\delta)$, since $\pi(x)$ will already be known.
(This is also true for the standard Metropolis method.)

Procedure~4 was shown in the previous section to leave $\pi$
invariant.  Accordingly, if we apply this procedure repeatedly, we are
justified in estimating expectations of functions of state using the
final states obtained at the end of each such application (after
discarding a suitable initial ``burn-in'' period).  This may be
inefficient, however.  Each application of Procedure~4 involves $K$
Metropolis updates.  If we did these in the standard manner, we would
be able to use all $K$ of these states when estimating expectations, not
just the last of them.  This will produce more accurate estimates,
although for difficult problems, in which these $K$ states are highly
dependent, the gain may be small (and may not be worth the cost of
storing these additional states).  

Fortunately, we can estimate expectations using all these states with
the short-cut Metropolis method as well.  In equilibrium, the state
before applying Procedure~4 will have distribution $\pi$.
Since every application of Step~(2) leaves $\pi$
invariant, the states at the end of Step~(2) will also have
distribution $\pi$, so we could use all $M$ of these states when
estimating expectations, not just the last.  Similarly, each
application of $T\met$ in Step~(2a) leaves $\pi$ invariant, so we can
decide to use all $K\!=\!ML$ of these states when estimating
expectations.  Note that if we decide to do this, we must use all
these states regardless of whether or not $T\met$ involved a
rejection, and whether or not the direction of simulation was reversed
in Step~(2a) --- any scheme in which the states to use are chosen
based on the results of the simulation might introduce bias.

If we decide to use all the $K$ states produced using $T\met$ (or if
we decide to use all $M$ states at the end of Step~(2)), we have
several options when one of these states turns out to be a duplicate
of an earlier state.  The simplest option is to simply copy the earlier
state to the area of memory reserved for the new state, or to write
the earlier state a second time to the output file, if the output is not
stored in main memory.  If states are stored in main memory, we could
instead just store a pointer to the earlier state.  Perhaps the most
efficient (though more complicated) method would be to store only a
single copy of each state, but to accompany this copy with a count of
how many times it should be included in the averages used to estimate
expectations.  (A similar issue arises with the standard Metropolis
method whenever a proposal is rejected, but since rejection rates are
usually not extreme, the usual practice of simply copying the rejected
state is generally adequate.)

If states are very large, and only the final states of short-cut
sequences are used for estimating expectations, it is possible to use
a method that avoids storing any but the initial and current states.
At the beginning of a short-cut sequence, the initial state and the
state of the pseudo-random number generator are saved.  Metropolis
updates are then simulated, with only the current state being
retained.  If a reversal occurs, the initial state may have to be
restored, at which point the state of the pseudo-random generator is
again saved.  During this procedure, we don't actually copy any states
when performing updates that produce previously-computed states (which
we haven't saved) --- we just keep track of where we are in the
earlier part of the sequence.  If the final state of the sequence is
not a copy of an earlier state, it will be available as the current
state.  Otherwise, we record the position of the final state that was
reached, save the final state of the pseudo-random number generator,
restore the initial state and the appropriate saved state of the
pseudo-random number generator, and simulate the number of Metropolis
updates needed to re-create the desired final state.  The final state
of the pseudo-random number generator is then restored.  This
procedure reduces the amount of memory required, but may increase the
computation time needed by up to a factor of two, though typically the
average increase will be much less than this.

Standard errors for estimates of expectations are usually found from
the sample variance and sample autocorrelation function (see, for
example, Neal 1993, Section 6.3).  If all $K$ states produced by a
short-cut Metropolis update are used, one might wonder whether this is
valid, since the process appears to be non-stationary, in which case
the autocorrelation function would not be well-defined.  In fact,
however, one can view the process as being stationary if the state is
extended to include auxiliary variables, including one that produces
alternation between Steps~(1) and~(2) of Procedure~4.  Alternatively,
one can look at the sample variance and sample autocorrelation
function of block averages, over the $K$ states in each short-cut
Metropolis update, or over longer blocks consisting of several
short-cut Metropolis updates with different stepsizes.

\section{\hspace*{-7pt}Demonstrations on simple
distributions}\label{sec-demo}\vspace*{-10pt}

To illustrate the operation of the short-cut Metropolis procedure, and
provide some insight into its performance, I will show how it works
when sampling from three simple distributions --- a one-dimensional
distribution for which the updates can easily be visualized, a
multivariate Gaussian distribution in which the variance in some
directions is much smaller than in others, and a ``funnel''
distribution that has features typical of Bayesian hierarchical
models.  

The program for these examples is listed in the Appendix, and
is available from my web page.

\subsection{\hspace*{-4pt}A one-dimensional mixture distribution}\vspace*{-4pt}

I will start with a one-dimensional example, since this allows for
easy visualization, though performance on one-dimensional examples is
not always typical of what happens in higher dimensions.

The distribution used here is an equal mixture of two Gaussian
distributions, one with mean 0 and standard deviation 10, the other
with mean 10 and standard deviation 1.  The density function for this
mixture is
\beq
   \pi(x) & = & {1 \over 2} {1 \over \sqrt{2\pi} \, 10} \exp(-(x/10)^2/2)
              \ +\ {1 \over 2} {1 \over \sqrt{2\pi}} \exp(-(x-10)^2/2)
\eeq
We will try to estimate the mean of this distribution from points
generated using Metropolis updates.  The true value of the mean is
$(1/2)0+(1/2)10 = 5$.

I used a Gaussian proposal distribution centred on the current state,
with standard deviation $w$ (in other words, $\delta$ was Gaussian
with mean zero and standard deviation 1, so that $w\delta$ had
standard deviation $w$).  If the distribution to sample from were
Gaussian with standard deviation~1, a stepsize of $w=2$ might be
appropriate, while if the distribution to sample from were Gaussian
with standard deviation 10, we might choose $w=20$.  Since the actual
distribution is a mixture of these two, we might be uncertain whether
the best stepsize is $w=2$ or $w=20$, or we might think that we need
to use both stepsizes at different times.

I tried five sampling methods on this distribution:\vspace*{-8pt}
\begin{quotation}\noindent
Standard Metropolis with $w=2$ \\[2pt]
Standard Metropolis with $w=20$ \\[2pt]
Naive adaptive Metropolis: \\
\hspace*{30pt}use $w=2$ if 
              there were more than 5 rejections in the last 10 updates \\
\hspace*{30pt}use $w=20$ otherwise \\[2pt]
Short-cut Metropolis with $l=0$ and $h=L-1$, alternating between two
sequences: \\
\hspace*{30pt}using $w=2$ with $L=5$ and $M=6$ (so $K=30$) \\
\hspace*{30pt}using $w=20$ with $L=5$ and $M=18$ (so $K=90$) \\[2pt]
Short-cut Metropolis with $l=1$ and $h=L-1$, alternating between 
two sequences: \\
\hspace*{30pt}using $w=2$ with $L=5$ and $M=12$ (so $K=60$) \\
\hspace*{30pt}using $w=20$ with $L=5$ and $M=12$ (so $K=60$)\vspace*{-8pt}
\end{quotation}
The standard and naive adaptive Metropolis methods were run for $1.2$ 
million iterations,
and therefore required $1.2$ million evaluations of $\pi(x)$.
The short-cut Metropolis method with $l=0$ (reversing only when all
updates in a group were rejected) was run for $16500$ pairs of sequences
with $w=2$ and $w=20$, which also required about $1.2$ million evaluations of
$\pi(x)$.  The short-cut Metropolis method with $l=1$ (reversing when either
all or none of the updates in a group were rejected) was run for $18000$
pairs of sequences with $w=2$ and $w=20$, again requiring about $1.2$ million
evaluations of $\pi(x)$.  All states generated by each method were
averaged to estimate the mean of $x$.  For the standard and naive Metropolis
methods, there were $1.2$ million states.  For the two short-cut Metropolis
methods, the numbers of states averaged was $1.98$ million and $2.16$ million,
but many of these states were copies of other states.

The results are shown in Table~\ref{tbl-1dres}.  The rejection rate
for Metropolis updates includes, for the short-cut methods, those that
are not actually performed.  The autocorrelation time for $x$ is
defined to be one plus twice the sum of the autocorrelations for $x$
at lags one to infinity; it is here estimated using the estimated
autocorrelations at lags one to 500.  The estimated mean for $x$ is
the sample average for all states.  The standard error for this
estimate is found from the variance of the mixture distribution, which
is known to be exactly 75.5, and the effective sample size, which is
the number of states used for the average divided by the autocorrelation time.
For further details on computation of MCMC standard errors, see 
(Neal 1993, Section~6.3).

\begin{table}[t]
\begin{center}\begin{tabular}{lccrcc}\hline
     & millions & rejection & \multicolumn{1}{c}{autocorrelation} & estimated 
& standard error \\[-3pt]
     & of states &   rate    & \multicolumn{1}{c}{time} &  mean    
&   for mean \\[1pt]
\hline\\[-10pt]
Standard, $w=2$  & 1.20 & 0.274 & 153.6~~~~~~ & 5.020 & 0.098 \\
Standard, $w=20$ & 1.20 & 0.699 &  10.2~~~~~~ & 4.989 & 0.025 \\
Naive adaptive   & 1.20 & 0.531 &  14.9~~~~~~ & 6.002 & 0.031 \\
Short-cut, $l=0$ & 1.98 & 0.590 &  53.0~~~~~~ & 4.923 & 0.045 \\
Short-cut, $l=1$ & 2.16 & 0.487 & 105.1~~~~~~ & 5.033 & 0.061 
\end{tabular}\end{center}

\caption[]{Results of five Metropolis methods on 
           the one-dimensional mixture distribution.}\label{tbl-1dres}
\end{table}

Although standard Metropolis with stepsize $w=2$ would be good for
exploring the mixture component with standard deviation one, we can
see from its high autocorrelation time (and consequent large standard
error) that it does not move around the whole mixture distribution
efficiently.  The low overall rejection rate of 0.274 is another
indication that $w=2$ is too small.  The results with $w=20$ are much
better.  For both methods, the estimated mean for $x$ differs from the
true value of 5 by an amount that is compatible with the estimated
standard error.

We might hope to improve on both of these standard Metropolis methods
using an adaptive scheme that chooses between $w=2$ and $w=20$
according to whether the number of rejections in the last 10 updates
was more than 5 or not.  This naive adaptive scheme does produce a
reasonable rejection rate, and a fairly low autocorrelation time.  But
the estimate it produces for the mean of $x$ is completely wrong.  This
estimate of $6.002$ differs from the true value of 5 by more than 30
times the standard error, showing that this method is heavily biased.

The two short-cut Metropolis methods do get the right answer, within
plus or minus twice the standard error.  They both are more efficient
than standard Metropolis with $w=2$, but they are less efficient than
standard Metropolis with $w=20$.  This is not too surprising, since for
very low-dimensional distributions, large stepsizes can be desirable,
even when they produce large rejection rates (often even when the
rejection rate is much larger than the fairly moderate value of 0.699
seen here for $w=20$).  We will have to look at higher-dimensional
problems to assess the advantages of short-cut Metropolis.

\begin{figure}[p]

\hspace{0.51in} $w=2$ \hspace{0.97in} $w=20$
\hspace{1.0in} $w=2$ \hspace{0.97in} $w=20$ \vspace*{-30pt}

\hspace*{-15pt}\includegraphics{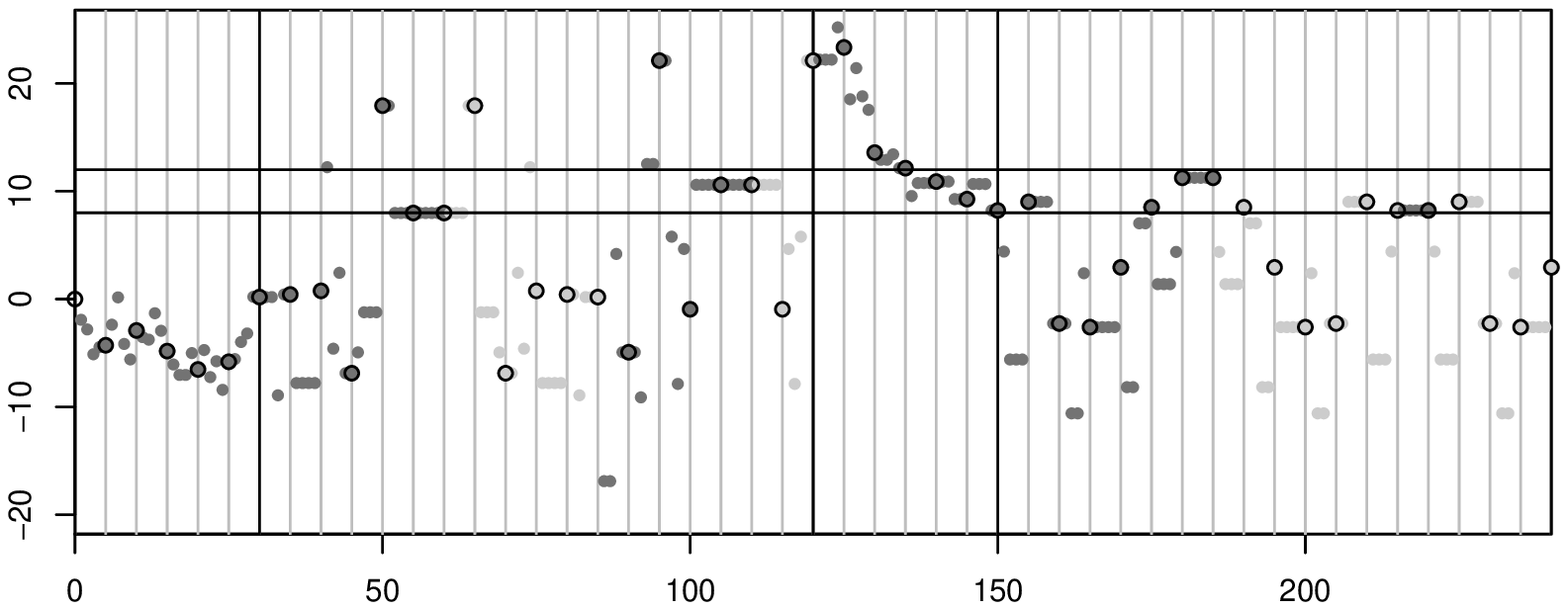}

\vspace*{-27pt}

\begin{center} \em 
  Two pairs of short-cut Metropolis sequences with $l=0$ and $h=L-1$
\end{center}

\vspace*{25pt}

\hspace{0.9in} $w=2$ \hspace{0.98in} $w=20$
\hspace{0.95in} $w=2$ \hspace{1.0in} $w=20$ \vspace*{-30pt}

\hspace*{-15pt}\includegraphics{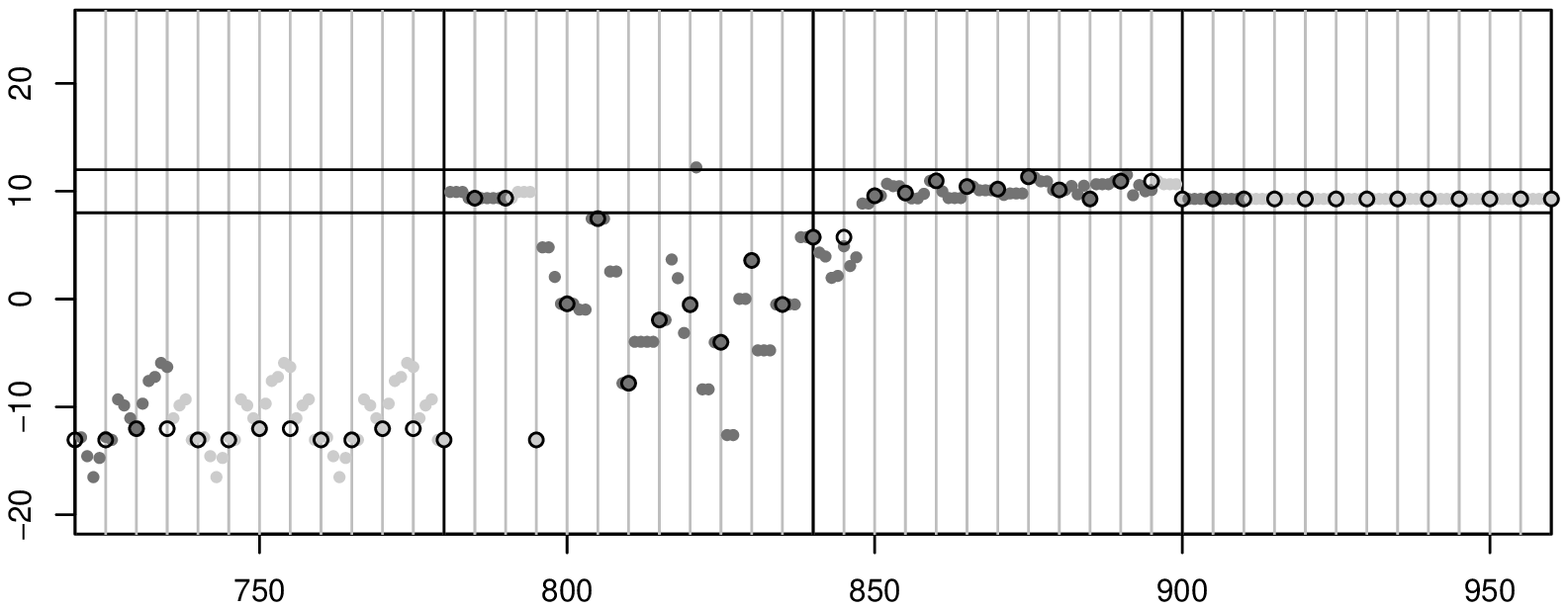}

\vspace*{-27pt}

\begin{center} \em 
  Two pairs of short-cut Metropolis sequences with $l=1$ and $h=L-1$
\end{center}\vspace*{10pt}

\caption[]{Short-cut Metropolis for the one-dimensional mixture distribution.
The top plot shows a portion of a short-cut Metropolis run
with $l=0$ and $h=L\!-\!1$, for which reversals
occur only when all updates in a group are rejected; the bottom plot
shows a portion of a run with $l=1$ and $h=L\!-\!1$,
for which reversals occur when either all or none of the updates in a group
are rejected.  The vertical axis is the state ($x$), 
with $10\pm2$ shown
by horizontal lines.  The horizontal axis indexes the updates, with the 
gray vertical lines marking the ends of groups, and the black vertical lines 
marking the ends of short-cut sequences.  Dark gray 
dots are states that required computation to find; light gray dots are states 
that were copied from earlier states.  Black circles are states at the ends
of groups.  The stepsizes used for the short-cut sequences are shown above
the plots.}\label{fig-1dpic}

\end{figure}

This one-dimensional example does allow for a good
visualization of short-cut Metropolis, as seen in
Figure~\ref{fig-1dpic}.  The top panel of this figure shows the first
four short-cut Metropolis sequences from the run in which reversals
were done only when a group of $L=5$ updates were all rejections.  In
the first sequence, with $K=30$ and $w=2$, all the states are outside
the narrow mixture component around $x=10$.  In this region, $w=2$ is
a very small stepsize.  Consequently, none of the groups consists only
of rejections, there are no reversals, and the entire sequence looks
exactly like a sequence of standard Metropolis updates.  The next
short-cut sequence, with $K=90$ and $w=20$, enters the region around
$x=10$, where there is a peak in the mixture density, which results in
a high probability of rejecting a proposal generated with $w=20$.  The
group of updates from index 56 to 60 are all rejected, resulting in a
reversal.  The states copied after that are shown as light dots.  Once
the initial state from index 30 has been copied, simulation of new
states (shown as darker dots) resumes, continuing until another group
of all rejections occurs at indexes 106 to 110, after which no further
computation of states is need to the end of this sequence.  The third
and fourth short-cut sequences repeat this pattern --- the first (with
$w=2$) has no reversals, and the second (with $w=20$) has two
reversals, with about half the states copied from other states.

The bottom panel in Figure~\ref{fig-1dpic} shows four short-cut
sequences from the run in which reversals occur when a group of
updates has either no rejections or all rejections.  In the first
sequence, with $w=2$ and $K=60$, the first group has no rejections.
The direction of simulation is therefore reversed, and the state is
restored to the initial state.  One update in the second group is a
rejection, but the third group has no rejections, causing a second
reversal. Note that at this reversal the state is restored to what it
was at the beginning of the third group.  Consequently, the state
after simulating this group is different from the state produced by
the last update in the group --- a situation that cannot happen when
reversals occur only when all updates in a group are rejected.
Subsequent groups in this sequence are just copies of
previously-computed states.  Effectively, the short-cut procedure has
``decided'' that the stepsize is too small, so that little work should
be done for this sequence.  In the second sequence, with $w=20$ and
$K=60$, a reversal occurs after the second group of updates, all of
which were rejected, but no second reversal occurs, so only a few
states were copied from other states.  Here, the short-cut method has
``decided'' that the stepsize is about right.  In the third sequence,
the stepsize is again $w=2$, but the chain enters the region around
$x=10$, for which this stepsize works well, so most of the states are
not copies.  The fourth sequence, with $w=20$, is also in the region
around $x=10$, for which $w=20$ is too large a stepsize.  The first
two groups consist only of rejections, so the entire sequence of
states is the same as the initial state, with only the first two
groups requiring any computation.

\subsection{\hspace*{-4pt}A multivariate Gaussian distribution}\vspace*{-4pt}

As a second demonstration, I applied standard and short-cut Metropolis
methods to the problem of sampling from a seven-dimensional
multivariate Gaussian distribution.  The mean of this distribution was
the zero vector, and its covariance matrix was diagonal, with the
variances of the first two components being 1 and the remaining five
being $0.1^2$.  The proposal distribution was also Gaussian, with mean
equal to the current point, and covariance matrix of $wI$.  Note that
since the proposal distribution is spherically symmetric, behaviour
would be unchanged if the distribution were rotated.  Accordingly,
although the seven coordinates are independent in the distribution
actually used, the same behaviour would be seen for any multivariate
Gaussian distribution whose covariance matrix has eigenvalues equal to
the variances used here.

I applied the standard Metropolis method to this distribution using
stepsizes of $w=0.02$, $w=0.1$, and $w=0.5$, in each case for 900000
updates.  The first of these stepsizes is too small, and the last is
too big, but we imagine that we do not realize this initially.
Accordingly, I also did a run in which these three stepsizes were
applied in turn, each for 200 updates at a time, with this cycle being
repeated 1500 times, again for a total of 900000 updates.  All 900000
states produced were used to estimate expectations.

I tried three versions of short-cut Metropolis.  All versions used
groups of $L=6$ updates, and cycled among sequences with $w=0.02$,
$w=0.1$, and $w=0.5$.  The total number of Metropolis updates for
these three stepsizes (including updates that were copied from
previous states) was always 600.  The number of times these three
sequences were repeated was adjusted so that the total number of
evaluations of $\pi(x)$ was approximately 900000, as for the standard
Metropolis runs.  More than 900000 states were produced (some copied
from earlier states), all of which were used to estimate expectations.

In the first version of short-cut Metropolis, reversals were done only
when all updates in a group were rejections (ie, $l=0$ and
$h=L\!-\!1$).  The short-cut sequences using $w=0.02$ were of length
$K=60$, those using $w=0.1$ were of length $K=150$, and those using
$w=0.5$ were of length $K=390$.  These increasing lengths were chosen
so that the method would spend most of its time using an appropriate
stepsize, regardless of which of the three was the appropriate one.
For sequences using the smallest stepsize of $w=0.02$, no reversals
were done, even if a group consisted only of rejections (ie, $l=0$ and
$h=L$), since there is no smaller stepsize to try in any case.  (In
other words, the 60 updates with $w=0.02$ were done in the standard
way, with no short-cut.)  The three sequences were repeated 4080
times, producing 2.448 million states.

The second version of short-cut Metropolis used sequences of length
$K=200$ for all three stepsizes.  Reversals were done when either all
or none of the updates in a group were rejections (ie, $l=1$ and
$h=L\!-\!1$), except that no reversal was done for a group of all
rejections when using the smallest stepsize, and no reversal was done
for a group of no rejections for the biggest stepsize.  For this version,
the three sequences were repeated 3000 times, producing 1.800 million states.

In the third version of short-cut Metropolis, reversals were done when
either all of the updates in a group were rejections, or when the
number of rejections was less than two (ie, $l=2$ and $h=L\!-\!1$) ---
except that, as before, reversals were done for the biggest stepsize
only for groups of all rejections, and for the smallest stepsize only
for groups with less than two rejections.  The rationale for this
version is that the optimal rejection rate is often somewhat greater
than 50\% (see Roberts and Rosenthal 2001), so using an asymmetrical
range for the desired number of rejections (from $l=2$ to $h=L-1$) may
be beneficial.  For this version, the three sequences were repeated
3720 times, producing 2.232 million states.

Results of estimating the expected value of the first component
of state (whose true value is zero) are shown in
Table~\ref{tbl-mvgres}.  Autocorrelation times were estimated using
the estimated autocorrelations up to lag 12000 for standard Metropolis
with $w=0.02$ and $w=0.5$, and to lag 8000 for the other methods.  The
standard errors shown account for the varying autocorrelation times,
and the fact that the short-cut Metropolis runs have a larger number
of states (found using the same number of evaluations of $\pi(x)$).
The actual differences between the estimates and the true mean of zero
are all consistent with the standard errors, exhibiting the usual
amount of chance variation.

\begin{table}[t]
\begin{center}\begin{tabular}{lccrcc}\hline
     & millions & rejection & \multicolumn{1}{c}{autocorrelation} & estimated 
& standard error \\[-3pt]
     & of states &   rate    & \multicolumn{1}{c}{time} &  mean    
&   for mean \\[1pt]
\hline\\[-10pt]
Standard, $w=0.02$    & 0.900 & 0.169 &  9677~~~~~~ & $+0.059$ & 0.104 \\
Standard, $w=0.1$     & 0.900 & 0.687 &  1271~~~~~~ & $+0.015$ & 0.038 \\
Standard, $w=0.5$     & 0.900 & 0.998 &  8311~~~~~~ & $-0.102$ & 0.096 \\
Standard, three $w$'s & 0.900 & 0.618 &  3998~~~~~~ & $-0.023$ & 0.067 \\
Short-cut, $l=0$      & 2.448 & 0.837 &  4719~~~~~~ & $+0.044$ & 0.044 \\
Short-cut, $l=1$      & 1.800 & 0.618 &  4427~~~~~~ & $-0.061$ & 0.050 \\
Short-cut, $l=2$      & 2.232 & 0.618 &  4729~~~~~~ & $+0.080$ & 0.046 
\end{tabular}\end{center}

\caption[]{Results of seven Metropolis methods on 
           the multivariate Gaussian distribution.}\label{tbl-mvgres}
\end{table}

The smallest standard error is for the standard Metropolis run with
$w=0.1$.  Standard Metropolis runs with $w=0.02$ and $w=0.1$ produced
much larger standard errors, showing that these stepsizes are not
suitable.  By assumption, however, we don't know ahead of time that
$w=0.1$ is the best stepsize.  If we therefore use all three stepsizes
in turn, we pay a price in terms of a larger standard error.  The cost
can be estimated by the square of the ratio of the standard errors,
which is $(0.067/0.038)^2 = 3.11$.  This measures how much longer a run
we would need using all three stepsizes to get the same accuracy as
when using just $w=0.1$.  As expected, it is near $3$, since only
$1/3$ of the time is spend using the good stepsize of $w=0.1$ in the
run using all three stepsizes.

We hope to do better than this using the short-cut method.  As seen in
Table~\ref{tbl-mvgres}, all three versions of the short-cut method
that were tried do indeed have smaller standard errors than standard
Metropolis using all three stepsizes, though their standard errors are
greater than that of standard Metropolis using $w=0.1$.  The
differences between the three short-methods are fairly small.  The
estimated advantages over the standard method using all three
stepsizes range from $(0.067/0.050)^2=1.80$ to $(0.067/0.044)^2=2.32$.

\begin{table}[t]
\begin{center}\begin{tabular}{rccc}\hline
                 & ~~$w=0.02$~~ & ~~$w=0.1$~~ & ~~$w=0.5$~~ \\[1pt]
\hline\\[-10pt]
            $K=$ &  60  &  150 & 390  \\[3pt]
Short-cut, $l=0$ & 0.00 & 0.09 & 0.95 \\[8pt]
            $K=$ & 200  &  200 & 200  \\[3pt]
Short-cut, $l=1$ & 0.49 & 0.13 & 0.90 \\
Short-cut, $l=2$ & 0.79 & 0.12 & 0.90 
\end{tabular}\end{center}

\caption[]{Fractions of states copied from earlier states for three
           versions of short-cut Metropolis.  
The lengths ($K$) of the
sequences for each stepsize ($w$) are shown as well.
}\label{tbl-cpy}
\end{table}

We can gain some insight into how the short-cut methods perform by
looking at the fraction of states that were copied from earlier states
(and hence required no evaluation of $\pi(x)$).  These fractions are
shown in Table~\ref{tbl-cpy}.  In the version in which reversals occur
only when a group consists of all rejections, almost all of the
updates with $w=0.5$ were copies, while almost none of the updates with
$w=0.1$ were copies.  (There were no copies for $w=0.02$, since these
updates were done with standard Metropolis.)  This is just what we
hoped for.  However, the updates done with $w=0.02$ still represent an
inefficiency, which is diminished but not eliminated by the fact
that the sequences with $w=0.02$ are shorter than those with $w=0.1$.

In the second version, reversals are done when either all or no
rejections occur in a group.  Again, most updates with $w=0.5$ are
copies.  Over half the updates with $w=0.02$ required actual
computation, however.  Since in this method the lengths of the
sequences for different stepsizes were the same, the result is that
the inefficiency from doing updates with $w=0.02$ was greater than for
the first version.

The third version reduces this problem by reversing when the number of
rejections in a group is either zero or one (except for the largest
stepsize of $w=0.5$).  Reversals then happen sooner with $w=0.02$, and
consequently more states are copies, and less time is wasted with this
unsuitable stepsize.  The standard error is therefore reduced.

\subsection{\hspace*{-4pt}A ``funnel'' distribution}\vspace*{-4pt}

As a final illustration, I will show how short-cut Metropolis can be
advantageous when no single stepsize is optimal for all regions of
the distribution being sampled.  I have previously used this example
to illustrate the advantages of slice sampling (Neal 2003).

The state for this example consists of ten real-valued components, $v$
and $x_1$ to $x_9$.  The marginal distribution of $v$ is Gaussian with
mean zero and standard deviation~3.  Conditional on a given value of
$v$, the variables $x_1$ to $x_9$ are independent, with the
conditional distribution for each being Gaussian with mean zero and
variance $e^v$.  The resulting distribution, $\pi$, has a shape
resembling a ten-dimensional funnel, with small values for $v$ at its
narrow end, and large values for $v$ at its wide end.  Such a
distribution is typical of priors for components of Bayesian
hierarchical models --- $x_1$ to $x_9$ might, for example, be random
effects for nine subjects, with $v$ being the log of the variance of
these random effects.  If the data happen to be largely uninformative,
the problem of sampling from the posterior will be similar to that of
sampling from the prior, so this test is relevant to actual Bayesian
inference problems.

I will focus here on estimating the mean of the $v$ component.  The
distribution is defined so that the true mean of $v$ is zero, but we
pretend here that we don't know that, and see how well various
Metropolis methods can estimate this mean, and more generally, the
marginal distribution of $v$.

The Metropolis methods I will consider update $v$ and $x_1$ to $x_9$
simultaneously, using a multivariate Gaussian proposal distribution
with mean equal to the current state and covariance matrix $wI$.  All
the methods employ sequences of 1000 Metropolis updates, with only the
last state from each sequence being used to estimate expectations.
The initial state had $v=0$ and all $x_i=1$.

\begin{figure}[p]

\vspace*{-10pt}

\hspace*{-5pt}\includegraphics{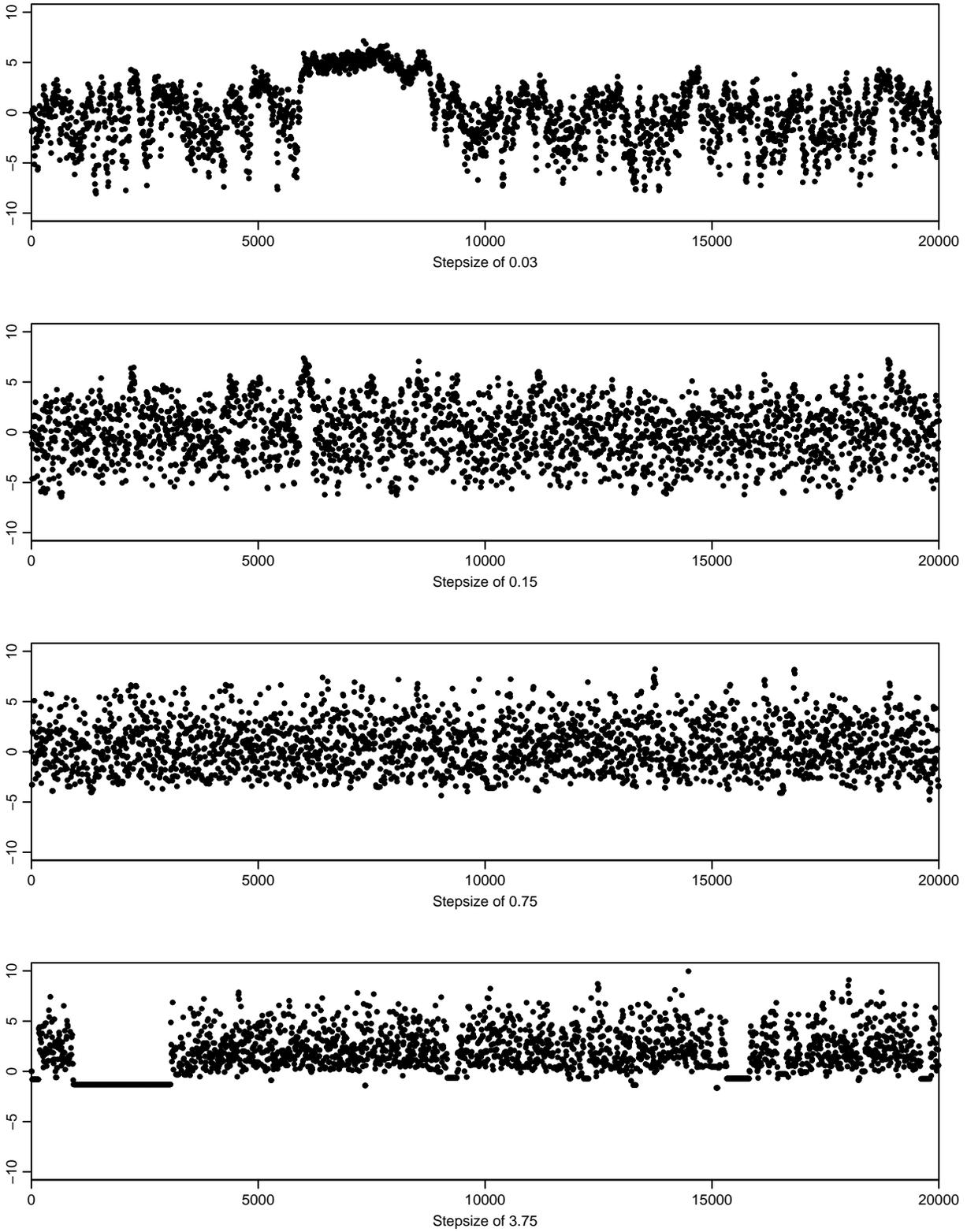}

\caption[]{The standard Metropolis method applied to the funnel distribution, 
using various stepsizes.  The horizontal axis indexes groups of 1000 standard 
Metropolis updates, with every eighth group being shown.  The vertical axis 
is the value of $v$ at the end of each group of updates.}\label{fig-stand3}

\end{figure}

Using standard Metropolis updates with a fixed value for
$w$ can be disastrous.  Figure~\ref{fig-stand3} shows values sampled
for $v$ using four stepsizes of $w=0.03$, $w=0.15$, $w=0.75$, and
$w=3.75$, in runs consisting of $20000$ sequences of 1000 updates.
When $w=3.75$, the run never goes much below zero.  This behaviour can
be explained by imagining what would happen if the chain reached a
point with a much smaller value of $v$, such as $v=-6$, along with
suitable values for $x_1$ to $x_9$.  The standard deviation of the
$x_i$ given this value of $v$ is $e^{-6/2} = 0.05$.  With $w=3.75$,
the probability of proposing a state with reasonable values of $x_1$
to $x_9$ is extremely small, so the probability of rejection is
extremely high, and the Markov chain will stay at the state with this
small value of $v$ for a very long time.  Since the chain leaves $\pi$
invariant, it follows that movement from a state with a large value of
$v$ to a state with a small value of $v$ must be very infrequent ---
so infrequent that, as we see in the figure, it never happens in 20
million Metropolis updates.  We do see signs of the problem in the
long strings of rejections when the chain enters states with $v$
slightly less than zero.

When $w=0.75$, the same problem exists.  None of the states sampled
have $v$ less than $-5$, which ought to occur with probability 0.048.
Worse, the plot in Figure~\ref{fig-stand3} for $w=0.75$ shows no sign
of there being any problem with his run.  This is the dreaded
``nightmare'' scenario of MCMC, in which lack of convergence leads to
drastically wrong results that appear to be correct.

Metropolis updates using the smallest stepsize of $w=0.03$ also do not
work well.  With this stepsize, small values of $v$ are well visited,
but large values, greater than about 4, are visited only in rare
extended excursions, such as the one in Figure~\ref{fig-stand3}
between sequences 6000 and 9000.  If only the first quarter of this
run had been done, no states with $v>5$ would have been sampled, with
no clear indication of a problem.  The run with $w=0.15$ appears much
better, sampling both states where $v$ is above 5 and states where it
is below $-5$.  A deficiency of very large and very small values for
$v$ is apparent, however.  In any case, it is hard to see how we would
know enough to choose $w=0.15$ ahead of time.

\begin{figure}[t]

\vspace*{-10pt}

\hspace*{-5pt}\includegraphics{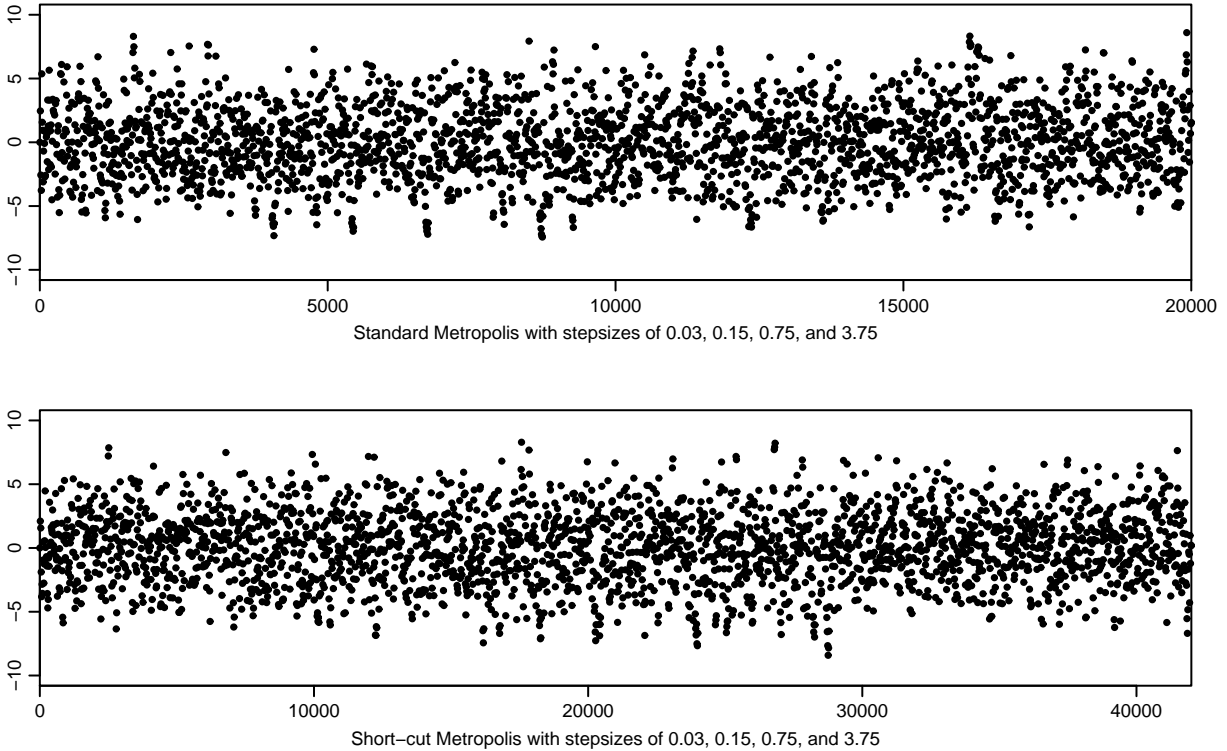}

\vspace*{-4.4in}

\caption[]{Standard and short-cut Metropolis methods using multiple
stepsizes, applied to the funnel distribution.  The plots are
analogous to those in Figure~\ref{fig-stand3}.  The run for the
short-cut method was about twice as long as the standard Metropolis
run (with the number of evaluations of $\pi$ being equal).
Accordingly, while the state after every eighth sequence is shown for
the standard method, every sixteenth is shown for the short-cut
method, so that the plots can be more easily compared.}\label{fig-ms3}

\end{figure}

Accordingly, we might try using standard Metropolis with a stepsize
that cycles amongst the values $w=0.03$, $w=0.15$, $w=0.75$, and
$w=3.75$.  The top plot of Figure~\ref{fig-ms3} shows a run of this
sort, in which these stepsizes are applied in turn for sequences of
1000 updates.  The results now appear to be acceptable, although there
is still some stickiness for very large or small values of
$v$.

We may hope to do better using the short-cut Metropolis method, by
avoiding the computation time spent on stepsizes that are not
appropriate.  A short-cut Metropolis run was done in which sequences
of 1000 Metropolis updates were done as $M=25$ groups of $L=40$
updates.  Reversals were done when a group consisted of all
rejections, and when a group had fewer than three rejections (ie,
$l=3$ and $h=39$), except that no reversals were done on all
rejections for the largest stepsize, or on fewer than three rejections
for the smallest stepsize.  Because some states in the short-cut run
are copied without further computation, 42000 sequences of 1000
updates can be done with the same number of evaluations of $\pi$ as
with the 20000 sequences done using standard Metropolis.

The bottom plot of Figure~\ref{fig-ms3} shows that the short-cut
method produces good results, sampling both large and small values of
$v$ at least as well as the standard Metropolis run using four
stepsizes.

Table~\ref{tbl-fun} provides a quantitative comparison of the methods.
The estimates in this table are based only on the final state from
each sequence of 1000 Metropolis updates, and autocorrelations are for
these final states, not for the states after each Metropolis update.
Autocorrelation times were estimated using estimated autocorrelations
up to lag 1000 for standard Metropolis with $w=0.03$, $w=0.15$, and
$w=3.75$, up to lag 100 for standard Metropolis with $w=0.15$, and up
to lag 50 for the standard and short-cut methods using all four stepsizes.

The standard Metropolis runs with $w=0.75$ and $w=3.75$ produce
estimates for the mean of $v$ that are far from the true value, much
further than would be expected given the computed standard errors.
These unrealistic standard errors result from underestimation of the
autocorrelation time, which for these runs is actually extremely
large.  The standard error for the standard Metropolis run with
$w=0.03$ may be realistic, but is uncomfortably large.

Standard Metropolis with $w=0.15$, standard Metropolis using all four
stepsizes, and short-cut Metropolis using these four stepsizes all
produce reasonable estimates, consistent with their standard errors.
The short-cut Metropolis method produced the best results.  Its
advantage over standard Metropolis using all four stepsizes is a factor
of $(0.090/0.073)^2 = 1.52$.  

\begin{figure}[b]

\vspace*{-5pt}

\vfill

\includegraphics{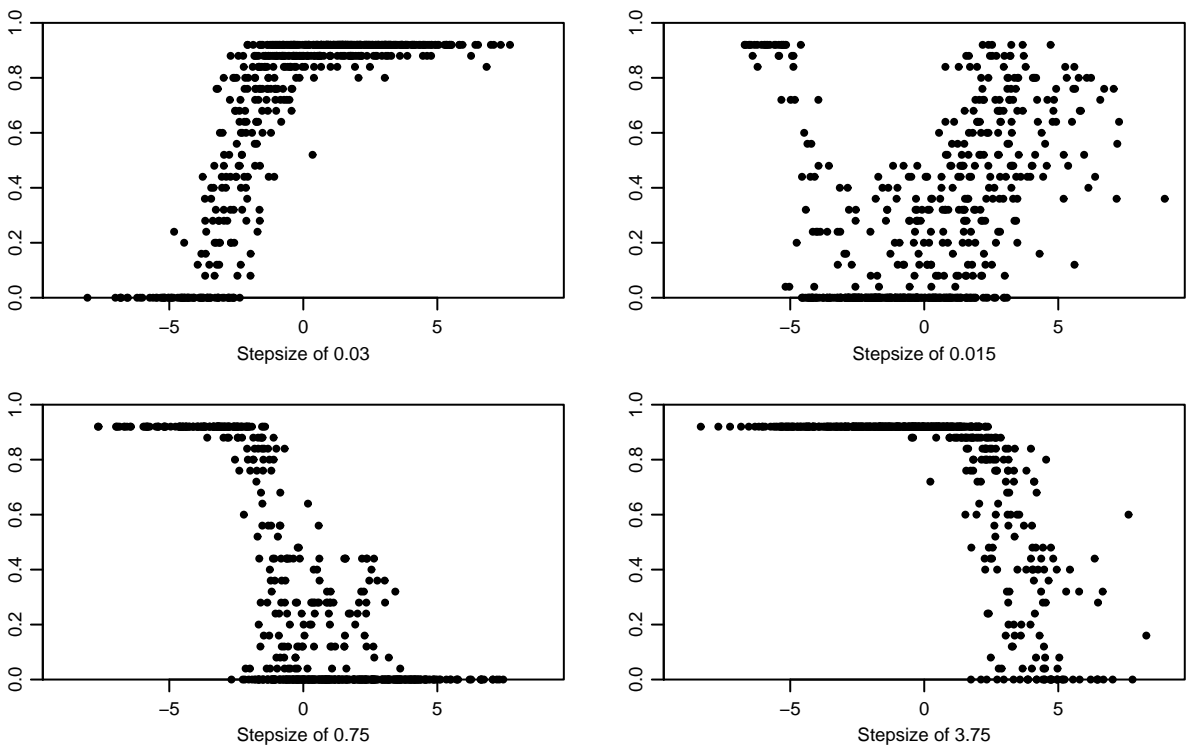}

\vspace*{-8pt}

\caption[]{How the fraction of states copied from earlier states
by short-cut Metropolis applied to the funnel distribution varies
with stepsize.
The horizontal axis is the value of $v$ at the start of a short-cut
sequence.  The vertical axis is the fraction of states in a short-cut
sequence that were copied from earlier states.  A point is plotted
for every fifth short-cut sequence in the run.\vspace*{-30pt}}\label{fig-cop3}

\end{figure}

Figure~\ref{fig-cop3} shows how this advantage was obtained.  For each
of the four stepsizes, a substantial fraction of states were sometimes
copied from earlier states, with the amount of copying varying with
the value of $v$ at the start of a sequence.  In effect, the short-cut
method allows the stepsize that is predominantly used to vary
depending on where in the state space the chain is currently located.

\begin{table}[t]
\begin{center}\begin{tabular}{lccrcc}\hline
     & thousands & rejection & \multicolumn{1}{c}{autocorrelation} & estimated 
& standard error \\[-3pt]
     & of states &   rate    & \multicolumn{1}{c}{time} &  mean    
&   for mean \\[1pt]
\hline\\[-10pt]
Standard, $w=0.03$    & 20 & 0.097 &  750~~~~~~ & $+0.143$ & 0.581 \\
Standard, $w=0.15$    & 20 & 0.324 &   39~~~~~~ & $+0.063$ & 0.133 \\
Standard, $w=0.75$    & 20 & 0.736 &    9~~~~~~ & $+0.501$ & 0.065 \\
Standard, $w=3.75$    & 20 & 0.968 &  438~~~~~~ & $+1.683$ & 0.444 \\
Standard, four $w$'s  & 20 & 0.540 &   18~~~~~~ & $+0.061$ & 0.090 \\
Short-cut, four $w$'s & 42 & 0.542 &   25~~~~~~ & $-0.022$ & 0.073 
\end{tabular}\end{center}

\caption[]{Results of several Metropolis methods applied to the funnel 
           distribution.}\label{tbl-fun}
\end{table}

\section{\hspace*{-7pt}Discussion}\label{sec-disc}\vspace*{-10pt}

In this paper, I have shown that by using the short-cut Metropolis
method the stepsize for the random-walk Metropolis method can in
effect be adaptively chosen from a small number of alternatives,
without the complications that arise in other adaptive schemes that
abandon the Markov property.  The short-cut method is also capable
of using different stepsizes in different regions of the state space.

Two general strategies have been demonstrated in the examples above.
In one strategy, short-cuts are taken only when the rejection rate is
high (due to the stepsize being too big).  We cyclicly simulate
short-cut sequences using a range of stepsizes, with the sequences
using bigger stepsizes being much longer than those using smaller
stepsizes.  If the biggest stepsize turns out to be best, it will
dominate the computation time simply because the sequence using it is
much longer than those using smaller stepsizes.  If instead a smaller
stepsize turns out to be best, it will dominate because the sequences
with bigger stepsizes will take little time to simulate, once two
groups of updates consisting of all rejections are encountered.  A
limitation of this strategy is that the sequences using increasing
stepsizes must have lengths that increase exponentially.  If we need
to use a large range of stepsizes, we might find that the length of
the longest sequence is greater than the total number of updates
we wish to perform.

The second strategy does not suffer from this problem.  It also
cyclicly simulates short-cut sequences using several stepsizes, but
these sequences can all be of the same length.  Short-cuts are taken
when the rejection rate is either too high (stepsize too big) or the
rejection rate is too low (stepsize too small).  If only one of the
stepsizes used is appropriate, only the sequence using this stepsize
will take appreciable time to simulate.  Simulations of sequences
using the other stepsizes will soon encounter two groups of updates
for which the number of rejections is either below the limit lower
established, or above the upper limit, after which no further
computation of states is required.  A fairly large number of stepsizes
can feasibly be used with this strategy.

With either strategy, it would be possible to pick stepsizes randomly
from some range, rather than cycling through a small number of
possible stepsizes.  The length ($K$) of the sequence could be a
function of the chosen stepsize, as could $l$, $h$, and $L$.

One possible improvement would be to look not at the actual number of
rejections within a group of updates, but rather at the expected
number of rejections, which is just the sum of one minus the
acceptance probability of equation~(\ref{eq-met-acc}) over all updates
in a group.  The set $A$ in equation~(\ref{eq-transA}) would be
defined in terms of the expected number of rejections for a group of
updates starting in the given state.  We must be careful, however.
For the method to remain valid, we must either look at whether the
state at the start of a group \textit{or} the state at the end of the
group is in $A$, or define the set $A$ so that if the start state is
in $A$ then the end state will be also.  The latter approach is
probably preferable.  We could, for example, define $A$ to consist of
states that produce a group of updates for which the \textit{average}
of the expected number of rejections when simulating the group
forwards and when simulating it backwards is outside some desired
interval.

This change would reduce the amount of random variation in whether or
not a group of updates causes a reversal in the simulation.  This is
likely to be beneficial, since it will lead to unsuitable stepsizes
being more reliably identified.  Situations where bad luck leads to
two early reversals even though the stepsize is suitable will be
less common.  Conversely, it will be less likely that a sequence using
an unsuitable stepsize will avoid encountering two reversals early on.

The examples in this paper use Metropolis updates that change all
components of the state at once.  It is common to instead perform a
sequence of Metropolis updates, each of which proposes to change only
one component (or sometimes a few components).  The simplest approach
to applying the short-cut technique in this context would use a single
stepsize parameter, $w$, multiplying the scale of all the proposal
distributions, and base reversals on the total number of rejections
for updates of all components.  This might, however, lead to most
computation time being devoted to a stepsize that is suitable for most
of the components, but is much too large or small for a one or a few
components.  One might instead base reversals on the maximum rejection
rate over all components, to avoid the possibility that one component
becomes ``stuck'' as a result of using too large a $w$.

Much better, however, would be to somehow adaptively choose different
stepsizes for different components.  Similarly, if we update all
components at once, using a Gaussian proposal distribution, we would
ideally adaptively choose the entire covariance matrix, not just a
single scale factor, $w$.  Unfortunately, the short-cut technique does
not easily handle a large number of tuning parameters.  We could
perform short-cut sequences using randomly-chosen values for these
parameters, and do reversals based on some criterion that indicates
whether the chosen values are unsuitable.  However, with many
parameters, suitable values might be chosen very rarely.  This appears
to be a fundamental limitation of the short-cut method.

Situations in which we wish to tune just one or a few parameters of an
MCMC method are not uncommon, however.  Aside from the Metropolis
methods discussed in this paper, one might try to extend the short-cut
technique to slice sampling (Neal 2003), combining slice sampling's
short-term adaptation within each update with the short-cut method's
longer-term ``adaptation''.  Exactly how to apply the short-cut method
to slice sampling remains to be worked out, however.  It is possible
to use the short-cut technique for MCMC methods based on simulation of
Hamiltonian dynamics (for a review, see Neal 1993, Section~5), where
the parameter to be tuned is the stepsize for a discretization of the
dynamics.  This stepsize must be kept small enough that the
discretization error in the Hamiltonian is not too large.  This is a
very natural setting for the short-cut method, since the dynamical
updates are already deterministic.  Indeed, it is in this context that
the short-cut idea first occurred to me.

An irony with adaptive schemes is that while they are designed to ease
the burden of setting the parameters of an MCMC method, they
themselves introduce even more parameters to be set.  Thus with the
short-cut method, we must select the set of $w$ values, the group
size, $L$, the number of groups, $M$, and the parameters, $l$ and $h$,
that define the desired range of rejections.  Of course, the hope is
that setting these parameters will be relatively easy, whereas picking
a suitable value for $w$ when we are ignorant about the
characteristics of the distribution $\pi$ may be impossible.  It will
take experience with a wide range of problems to see how well this
hope is fulfilled, and to develop guidelines for how best to use the
short-cut method in practice.

\newpage

\section*{Appendix --- An R function implementing short-cut 
          Metropolis}\vspace{-10pt}

The function below, written in R (see
\texttt{http://www.r-project.org}), implements a short-cut Metropolis
update.  This program (with additional code for gathering statistics)
is available, along with scripts for the demonstrations in this paper,
at \texttt{http://www.cs.utoronto.ca/$\sim$radford}.  Note that this
program is intended for demonstration purposes only.  I have made no
serious attempt to optimize the code.  Due to the interpretive
implementation of R, the short-cut method suffers from significant
overheads, which are not due to any fundamental aspect of the method.

The code below begins with a function for performing a standard Metropolis
update, followed by the function for short-cut Metropolis.  The version
here returns all $K$ states produced in the course of the short-cut
procedure, all of which may be used when estimating expectations of state.

\vspace*{10pt}

{\small

\begin{verbatim}
# DO ONE METROPOLIS UPDATE.
#
# Arguments:
#
#    initial.x    The initial state (a vector)
#    lpr          Function returning the log probability of a state, plus an
#                 arbitrary constant
#    pf           Function returning the random offset (a vector) for a proposal
#    w            The stepsize for proposals, multiplies the offset
#    initial.lpr  The value of lpr(initial.x).
#
# The value returned is a list containing the following elements:
#
#    next.x       The new state
#    next.lpr     The value of lpr(next.x)
#    rejected     TRUE if the proposal was rejected

metropolis.update <- function (initial.x, lpr, pf, w, initial.lpr)
{
  # Propose a candidate state, and evalute its log probability.

  proposed.x <- initial.x + w*pf()
  proposed.lpr <- lpr(proposed.x)

  # Decide whether to accept or reject the proposed state as the new state.

  if (runif(1)<exp(proposed.lpr-initial.lpr)) # accept
  { next.x <- proposed.x
    next.lpr <- proposed.lpr
    rejected <- FALSE
  }
  else # reject
  { next.x <- initial.x
    next.lpr <- initial.lpr
    rejected <- TRUE
  }

  # Return the new state, its log probability, and whether a rejection occurred.

  list (next.x=next.x, next.lpr=next.lpr, rejected=rejected)
}
\end{verbatim}

}

\newpage

{\small

\begin{verbatim}
# THE SHORT-CUT METROPOLIS METHOD.  Simulates a short-cut Metropolis update
# consisting of M*L Metropolis updates.
#
# Arguments:
#
#    initial.x    The initial state (a vector of length n)
#    lpr          Function returning the log probability of a state, plus an
#                 arbitrary constant
#    pf           Function returning the random offset (a vector) for a proposal
#    w            The stepsize for proposals, multiplies the offset
#    L            Number of updates in each group
#    M            Number of groups to simulate
#    min.rej      Minimum number of rejections for a good group of L updates 
#    max.rej      Maximum number of rejections for a good group of L updates
#
# The value returned is a list containing the following elements:
#
#    states1      A M*L+1 by n matrix whose rows contain the initial state 
#                 and the states after each Metropolis update, including
#                 updates in groups after which a reversal occurred
#    states2      A M+1 by n matrix whose rows contain the initial state and
#                 the states after each group of L Metropolis updates; the
#                 last row contains the final state from the whole sequence

short.cut.metropolis <- function (initial.x, lpr, pf, w, L, M, 
                                  min.rej=0, max.rej=L-1)
{
  # The function below puts together the list of results that we return.

  results <- function () list (states1=states1, states2=states2)

  # The function below performs the L Metropolis updates in a group, storing the
  # results in states1.  The last.x and last.lpr variables should contain the
  # previous state and its log probability; they are updated by this function.
  # The variable k indexes where in states1 the new states should be stored.  
  # It is incremented in this function.  The n.rejected variable is set to the 
  # number of the L updates that were rejections.  

  do.group <- function ()
  { n.rejected <<- 0
    for (l in 1:L)
    { update <- metropolis.update (last.x, lpr, pf, w, last.lpr)
      states1[k+1,] <<- last.x <<- update$next.x
      last.lpr <<- update$next.lpr
      k <<- k + 1
      n.rejected <<- n.rejected + update$rejected
    }
  }

  # Allocate space for the results.

  K <- M*L
  states1 <- matrix(NA,K+1,length(initial.x))
  states2 <- matrix(NA,M+1,length(initial.x))
\end{verbatim}\newpage\begin{verbatim}
  # Store the initial state in the first row of states1 and of states2, and 
  # evaluate its log probability.

  states1[1,] <- initial.x
  states2[1,] <- initial.x
  initial.lpr <- lpr(initial.x)
 
  # Do groups of Metropolis updates starting from the initial state, until we've
  # done many as were asked for, or until the number of rejections in a group 
  # is outside the limits.  States after each Metropolis update are saved
  # in states1.  The final state for each group is saved in states2.  Note 
  # that for a group with the number of rejections outside the limits, the 
  # state saved in states2 is not state the last state in states1, but rather
  # state from the start of the group.  The indexes in states1 of the start 
  # (upper0) and end (upper1) of these groups in are recorded for later use 
  # in copying states.

  last.x <- initial.x
  last.lpr <- initial.lpr
  k <- 1

  upper0 <- k
  while (k<=K)
  { states2[2+(k-1)/L,] <- last.x
    do.group()
    if (n.rejected<min.rej || n.rejected>max.rej)
    { upper1 <- k
      break
    }
    else
    { states2[1+(k-1)/L,] <- states1[k,]
    }
  }

  if (k>K) return (results())  # Return if we've done all M groups.

  # Copy already-computed states as "new" states as we move backwards through
  # the groups previously simulated.  Don't copy the last group causing the
  # reversal, for which the number of rejections was outside the limits.

  j <- upper1 - L - 1
  while (k<=K && j>=upper0)
  { states1[(k+1):(k+L),] <- states1[j:(j-L+1),]
    states2[1+(k+L-1)/L,] <- states2[1+(j-L)/L,]
    k <- k + L;  j <- j - L
  }

  if (k>K) return (results())  # Return if we've done all M groups.

  # Restore the initial state, then do more groups of Metropolis updates, 
  # until we've done as many as were asked for, or the number of rejections in
  # a group is outside the limits.  Record the indexes of the start (lower0) 
  # and end (lower1) of these groups in states1. 

  last.x <- initial.x
  last.lpr <- initial.lpr
\end{verbatim}\newpage\begin{verbatim}
  lower0 <- k
  while (k<=K)
  { states2[2+(k-1)/L,] <- last.x
    do.group()
    if (n.rejected<min.rej || n.rejected>max.rej)
    { lower1 <- k
      break
    }
    else
    { states2[1+(k-1)/L,] <- states1[k,]
    }
  }

  if (k>K) return (results())  # Return if we've done all M groups.

  # Copy already-computed states as "new" states, going back and forth over
  # the "lower" groups and the "upper" groups.

  repeat
  {
    # Copy the lower states backwards, excluding the group causing the reversal.

    j <- lower1 - L - 1
    while (k<=K && j>=lower0)
    { states1[(k+1):(k+L),] <- states1[j:(j-L+1),]
      states2[1+(k+L-1)/L,] <- states2[1+(j-L)/L,]
      k <- k + L;  j <- j - L
    }

    if (k>K) return (results())  # Return if we've done all M groups.

    # Copy the upper states forwards, including the group causing the reversal.

    j <- upper0+1
    while (k<=K && j<=upper1)
    { states1[(k+1):(k+L),] <- states1[j:(j+L-1),]
      states2[1+(k+L-1)/L,] <- states2[1+(j+L-2)/L,]
      k <- k + L;  j <- j + L
    }

    if (k>K) return (results())  # Return if we've done all M groups.

    # Copy the upper states backwards, excluding the group causing the reversal.

    j <- upper1 - L - 1
    while (k<=K && j>=upper0)
    { states1[(k+1):(k+L),] <- states1[j:(j-L+1),]
      states2[1+(k+L-1)/L,] <- states2[1+(j-L)/L,]
      k <- k + L;  j <- j - L
    }

    if (k>K) return (results())  # Return if we've done all M groups.
\end{verbatim}\newpage\begin{verbatim}   
    # Copy the lower states forwards, including the group causing the reversal.

    j <- lower0 + 1
    while (k<=K && j<=lower1)
    { states1[(k+1):(k+L),] <- states1[j:(j+L-1),]
      states2[1+(k+L-1)/L,] <- states2[1+(j+L-2)/L,]
      k <- k + L;  j <- j + L
    }

    if (k>K) return (results())  # Return if we've done all M groups.
  }
}
\end{verbatim}

}


\section*{Acknowledgements}\vspace{-10pt}

I thank David MacKay for helpful comments on the manuscript.  This
research was supported by the Natural Sciences and Engineering
Research Council of Canada.  I hold a Canada Research Chair in
Statistics and Machine Learning.

\section*{References}\vspace{-10pt}

\leftmargini 0.2in
\labelsep 0in

\begin{description}
\itemsep 2pt

\item
  Andrieu, C.\ and Moulines, E.\ (2005) ``On the ergodicity properties
  of some adaptive MCMC algorithms'', preprint obtainable from
  \texttt{http://www.maths.bris.ac.uk/$\sim$maxca/}.

\item
  Atchad\'{e}, Y.~F.\ and Rosenthal, J.\ S.\ (2005) ``On adaptive Markov
  chain Monte Carlo algorithms'', preprint obtainable from
  \texttt{http://www.mathstat.uottawa.ca/$\sim$yatch436/}.

\item
  Green, P.~J.\ and Mira, A.\ (2001) ``Delayed rejection in reversible
  jump Metropolis-Hastings'', {\em Biometrika}, Biometrika, vol.~88.


\item
  Haario, H., Saksman, E., and Tamminen, J.\ (2001) ``An adaptive Metropolis
  algorithm'', {\em Bernoulli}, vol.~7, pp.~223-242.

\item
  Liu, J.~S.\ (2001) \textit{Monte Carlo Strategies in Scientific Computing},
  Springer-Verlag.

\item
  Metropolis, N., Rosenbluth, A.~W., Rosenbluth, M.~N., Teller, A.~H., 
  and Teller, E.\ (1953) ``Equation of state calculations by fast computing 
  machines'', {\em Journal of Chemical Physics}, vol.~21, pp.~1087-1092.

\item
  Neal, R.~M.\ (1993) {\em Probabilistic Inference Using Markov Chain
  Monte Carlo Methods}, Technical Report CRG-TR-93-1, Dept.\
  of Computer Science, University of Toronto, 140 pages.  
  Obtainable from \texttt{http://www.cs.utoronto.ca/$\sim$radford/}.

\item
  Neal, R. M. (2003) ``Slice sampling'' (with discussion), 
  {\em Annals of Statistics}, vol.~1, pp.~705-767.

\item
  Roberts, G.~O.\ and Rosenthal, J.~S.\ (2001) ``Optimal scaling for 
  various Metropols-Hastings algorithms'', \textit{Statistical Science},
  vol.~16, pp.~351-367.

\item
  Tierney, L.\ (1994) ``Markov chains for exploring posterior distributions''
  (with discussion), \textit{Annals of Statistics}, vol.~22, pp.~1701-1762.

\item
  Tierney, L.\ and Mira, A.\ (1999) ``Some adaptive Monte Carlo methods 
  for Bayesian inference'', {\em Statistics in Medicine}, \textbf{18}, 
  2507-2515.

\end{description}

\end{document}